\documentclass[12pt]{article}
\usepackage{amssymb,latexsym,amsmath}
\usepackage{graphics}                 
\usepackage{color}                    
\usepackage{hyperref}                 
\usepackage[all]{xy}

\begin{document}
\pagestyle{plain} \headheight=5mm \topmargin=-5mm 

\title{Lawson homology, morphic cohomology and Chow motives}
\author{ Wenchuan Hu and Li Li}

\maketitle

\newtheorem{definition}{Definition}[section]
\newtheorem{theorem}{Theorem}[section]
\newtheorem{proposition}{Proposition}[section]
\newtheorem{lemma}{Lemma}[section]
\newtheorem{corollary}{Corollary}[section]
\newtheorem{question}{Question}
\newtheorem{remark}{Remark}[section]
\newtheorem{example}{Example}[section]

\def\s{\section}
\def\ss{\subsection}

\def\nn{\nonumber}
\def\bp{{\bf Proof.}\hspace{2mm}}
\def\qe{\hfill$\Box$}
\def\lj{\langle}
\def\rj{\rangle}
\def\dd{\diamond}
\def\ox{\mbox{}}
\def\lb{\label}
\def\rel{\;{\rm rel.}\;}
\def\vp{\varepsilon}
\def\ep{\epsilon}
\def\mod{\;{\rm mod}\;}

\def\dim{{\rm dim}}
\def\im{{\rm im}\;}

\def\span{{\rm span}}
\def\Griff{{\rm Griff}}
\def\sign{{\rm sign}\;}
\def\Supp{{\rm Supp}\;}
\def\Sp{{\rm Sp}\;}
\def\ind{{\rm ind}\;}
\def\rank{{\rm rank}\;}

\def\Na{{\cal N}}
\def\det{{\rm det}\;}
\def\dist{{\rm dist}}
\def\deg{{\rm deg}}
\def\Tr{{\rm Tr}\;}
\def\ker{{\rm ker}\;}
\def\Vect{{\rm Vect}}
\def\H{{\bf H}}
\def\L{{\mathbb{L}}}
\def\K{{\rm K}}
\def\R{{\mathbb{R}}}
\def\C{{\mathbb{C}}}
\def\P{{\rm P}}
\def\Q{{\mathbb{Q}}}
\def\Z{{\mathbb{Z}}}
\def\N{{\bf N}}
\def\F{{\bf F}}

\def\x{{\bf x}}
\def\y{{\bf y}}

\def\E{{\rm E}}
\def\J{{\cal J}}
\def\ave{{\textrm{ave}}}
\def\Ch{{\rm Ch}}

\begin{abstract}
In this paper, the Lawson homology and morphic cohomology are
defined on the Chow motives. We also define the rational coefficient 
Lawson homology and morphic
cohomology of the Chow motives of finite quotient projective
varieties. As a consequence, we obtain a formula for the Hilbert
scheme of points on a smooth complex projective surface. Further
discussion concerning generic finite maps is given. As a result, we
give examples of self-product of smooth projective curves with
nontrivial Griffiths groups by using a result of Ceresa.

\end{abstract}

\tableofcontents

\section{Introduction}\label{sec:introduction}

The main purpose in this paper is to define the Lawson homology and
morphic cohomology on the usual Chow motives as well as the Chow
motives of finite quotient projective varieties.

The Lawson homology groups (resp. morphic cohomology groups) are the
homotopy groups of the space of algebraic cycles (resp. algebraic
cocycles), first studied by Friedlander and Lawson. We explain
briefly their idea:

Let $X$ be a complex projective variety and ${\cal Z}_p(X)$ be the
abelian group of algebraic cycles of dimension $p$ on $X$. There is
a natural topology, namely Chow topology, on this abelian group
which is independent of the projective embedding of $X$. The Lawson
homology $L_pH_k(X)$ is defined to be the homotopy group
$$L_pH_k(X) :=\left\{
                \begin{array}{ll}
                  \pi_{k-2p}({\cal Z}_{p}(X)), & \hbox{if $ k\geq 2p$} \\
                  0, & \hbox{if $k<2p$}
                \end{array}
              \right.
 $$
 (cf.
\cite{Friedlander1}, \cite{Lawson1}, \cite{Lawson2}). The
topological group $\mathcal{Z}^p(X)$ of all algebraic cocycles of
codimension-$p$ on $X$ is defined as a homotopy quotient completion
(cf. \cite{FL1}, Definition 2.8)
$$\mathcal{Z}^q(X):=[\mathfrak{Mor}(X,\mathcal{C}_0(\mathbb{P}^q))/
\mathfrak{Mor}(X,\mathcal{C}_0(\mathbb{P}^{q-1}))]^+=\mathfrak{Mor}(X,\mathcal{Z}_0(\mathbb{A}^q)).$$
Take the $(l-2q)$-th homotopy group of the space of algebraic
cocycles instead of algebraic cycles, we get the morphic cohomology
$L^qH^l(X)$. (Partial motivation to study the homotopy of the cycle
space is Almgren's isomorphism, which asserts that for a topological
space $X$ satisfying reasonable conditions,
$$H_k(X)\cong\pi_{k-r}(Z_r(X)),$$
where $Z_r(X)$ is the space of $r$-dimensional integral cycles,
i.e., integral currents without boundary.)

Now we fix our notation of Chow motives. Let $\mathcal{V}$ denote
the category of (not necessarily connected) complex smooth
projective varieties. Given two smooth projective varieties $X$ and
$Y$. Suppose $X=\coprod X_\alpha$ is the decomposition of $X$ into
irreducible components. The group of correspondences of degree $r$
from $X$ to $Y$ is defined as
$$Corr^r(X,Y):=\oplus\, Ch^{\dim X_\alpha+r}(X_\alpha\times Y),$$
moreover, its tensor with $\mathbb{Q}$ is denoted by
$Corr^r_\mathbb{Q}(X,Y)$. The composition of two correspondences
$f\in Corr^r(X,Y)$ and $g\in Corr^s(Y, Z)$ gives a correspondence in
$Corr^{r+s}(X, Z)$. A correspondence $\mathbf{p}\in Corr^0(X,X)$ is
called a projector of $X$ if $\mathbf{p}^2=\mathbf{p}$. The category
of Chow motives $CH\mathcal{M}$ is given as follows (cf. \cite{CH}
for the version we used here): Objects in $CH\mathcal{M}$ are
triples $(X,\mathbf{p},r)$, or denoted by $h(X,\mathbf{p})(-r)$,
where $X\in\mathcal{V}$ , $\mathbf{p}$ is a projector of $X$,
$r\in\mathbb{Z}$. In particular, the motive $h(X, id_X)(-r)$ is simply denoted by
$h(X)(-r)$.  Morphisms are defined as
$$Hom_{CH\mathcal{M}}\big{(}(X,\mathbf{p},r), (Y,\mathbf{q},s)\big{)}:=\mathbf{q}\circ Corr^{s-r}(X, Y)\circ \mathbf{p}.$$
The composition of morphisms is defined as the composition of
correspondences.

The following result states a relation of motives and the morphic
cohomology. (Analogous result holds for Lawson homology, cf. Theorem
\ref{Th4.3} (i).)

\begin{theorem}[Theorem \ref{main1}]\footnote{This theorem is implicit, in a different formulation, in
\cite{Nenashev-Zainoulline},  since morphic cohomology gives an
example of oriented cohomology theory. But for the completeness of
the paper and the convenience of the readers, we still state it
here.} Given any $q, l\in \mathbb{Z}$, the morphic cohomology
$L^qH^l$ defines a covariant functor from the category
$\mathcal{CHM}$ to the category of abelian groups as follows:
$$L^qH^l(X,\mathbf{p},r):=\mathbf{p}_*(L^{q+r}H^{l+2r}(X)) \subseteq
L^{q+r}H^{l+2r}(X).$$ Given a morphism $\Gamma\in
Hom_{CH\mathcal{M}}\big{(}(X,\mathbf{p},r),
(Y,\mathbf{q},s)\big{)}$, the morphism
$$L^qH^l(\Gamma): L^qH^l(X,\mathbf{p},r)\to L^qH^l(Y,\mathbf{q},s)$$ is defined as the
restriction of the map
$$\Gamma_*: L^{q+r}H^{l+2r}(X)\to L^{q+s}H^{l+2s}(Y).$$
\end{theorem}

The advantage of this theorem is that we can apply results on
motives to morphic cohomology theory (and Lawson homology). Examples
are: the projective bundle theorem (Corollary \ref{projbundlethm})
which is firstly proved by Friedlander and Gabber in \cite{FG} , and
the blowup formula for Lawson homology (Corollary \ref{cor4.2}),
which is proved by the first author in \cite{author1}, a result of
Lima-Filho (\cite{Lima3}) for projective manifolds
 admitting cell-decompositions. By applying the above theorem to a
result of N. A. Karpenko in \cite{Karpenko}, we get certain decomposition
of Lawson homology and the morphic cohomology (Corollary \ref{cor4.3}).

\medskip
Then we discuss finite quotients of smooth (quasi-)projective
varieties.

Our first observation is a natural relation between the rational
coefficient Lawson homology of a smooth quasi-projective variety and
the one of its quotient.
\begin{proposition}[Proposition \ref{G-invariant}]
Let $\pi:X\to X':=X/G$ denote the quotient
 map of a quasi-projective variety  with a faithful action of a finite group $G$. Then there is a
canonical isomorphism
$$\pi_*: (L_pH_k(X,\mathbb{Q}))^G\cong L_pH_k(X',\mathbb{Q}), \quad \hbox{ for any $p, k\in \mathbb{Z}$}.$$
 and an isomorphism, when $X$ is projective, as follows
$$\pi_!: (L^qH^l(X,\mathbb{Q}))^G\cong L^qH^l(X',\mathbb{Q}), \quad \hbox{ for any $q, l\in \mathbb{Z}$}.$$
\end{proposition}

\begin{remark}
 Friedlander and Walker proved the proposition (in the proof of Theorem
 5.5 in \cite{FW}) under the assumption of the smoothness of the quotient $X'=X/G$.
\end{remark}

The following is our main result for quotient varieties of smooth
projective varieties by a finite group action. (Analogous result
holds for Lawson homology, as stated in Theorem \ref{Th4.3} (ii).)
\begin{theorem}[Theorem \ref{main2}]
  Given any $q, l\in \mathbb{Z}$, the $\mathbb{Q}$-coefficient morphic cohomology $L^qH^l(-,\mathbb{Q})$
defines a covariant functor from the category $CH\mathcal{M}'$ of
Chow motives of quotient varieties to the category of abelian groups
as follows:
$$L^qH^l((X',\mathbf{p},r),\mathbb{Q}):=\mathbf{p}_*(L^{q+r}H^{l+2r}(X',\mathbb{Q})) \subseteq
L^{q+r}H^{l+2r}(X',\mathbb{Q}).$$ Given a morphism $\Gamma\in
Hom_{CH\mathcal{M}'}\big{(}(X',\mathbf{p},r),
(Y',\mathbf{q},s)\big{)}$, the morphism
$$L^qH^l(\Gamma,\mathbb{Q}): L^qH^l((X',\mathbf{p},r),\mathbb{Q})\to L^qH^l((Y',\mathbf{q},s),\mathbb{Q})$$ is defined as the
restriction of the map
$$\Gamma_*: L^{q+r}H^{l+2r}(X',\mathbb{Q})\to L^{q+s}H^{l+2s}(Y',\mathbb{Q}).$$
\end{theorem}

As an application, we give a decomposition of the Lawson homology
and morphic cohomology for the Hilbert scheme $X^{[n]}$ of $n$
points on a smooth complex projective surface $X$. It is well-known
that $X^{[n]}$ is nonsingular (cf. \cite{Fogarty}). Let $X^{(n)}$ be
the $n$-th symmetric product of $X$ and let $\pi:X^{[n]}\rightarrow
X^{(n)}$ be the natural morphism, namely the Hilbert-Chow morphism.  We denote by
$\mathfrak{P}(n)$ the set of partitions of $n$. Any $\nu\in
\mathfrak{P}(n)$ determined a quotient variety $X^{(\nu)}$ which is
a product of symmetric products of $X$ (for detailed meaning of
notations appeared here, the reader is referred to \S
\ref{sec:Hilbert} or \cite{dCM}). Apply the above theorem to the
motivic decomposition of $X^{[n]}$ proved by de Cataldo and
Migliorini in \cite{dCM}, we obtain the following theorem:

\begin{theorem}[Theorem \ref{main3}]
Let $X$ be a smooth complex projective surface. Then there is an
isomorphism of Lawson homology groups for all $p,k\in \mathbb{Z}$:
$$\bigoplus_{\nu\in
\mathfrak{P}(n)}L_{p-n+l(\nu)}H_{k-2n+2l(\nu)}(X^{(\nu)},
\mathbb{Q})\longrightarrow L_pH_k(X^{[n]},\mathbb{Q})
$$
and an isomorphism of morphic cohomology groups
$$
\bigoplus_{\nu\in
\mathfrak{P}(n)}L^{q-n+l(\nu)}H^{l-2n+2l(\nu)-d_{\nu}}(X^{(\nu)},
\mathbb{Q})\longrightarrow L^qH^l(X^{[n]},\mathbb{Q}).
$$
\end{theorem}

\medskip

To state the further consequences, we need to introduce some
notations. The continuous homomorphism ${\cal Z}_p(X)\hookrightarrow
Z_{2p}(X)$  induces the cycle class map
$$ \Phi_{p,k}:L_pH_{k}(X)\rightarrow H_{k}(X). $$
Define $L_pH_{k}(X)_{hom}:={\rm ker\, } \Phi_{p,k}$. The special
case $L_pH_{2p}(X)_{hom}$ gives the {Griffiths group} of $p$-cycles
${\rm Griff}_p(X):={\mathcal Z}_p(X)_{hom}/{\mathcal Z}_p(X)_{alg}$.
(The isomorphism $L_pH_{2p}(X)_{hom}\cong {\rm Griff}_p(X)$ is shown
by Friedlander \cite{Friedlander1}.)

Then we give the following applications using the idea of motives
and correspondences: a new proof of a result of the first author
that the Lawson homology groups $L_1H_k(-)_{hom}$ and
$L_{n-2}H_k(-)_{hom}$ are birational invariants; Some properties of
the Lawson homology groups of unirational threefolds and fourfolds;
\emph{examples of self-products of generic curves carrying
nontrivial Griffiths groups}.

\begin{proposition}[Proposition \ref{nonzeroGriff}]
Let $C$ be generic smooth projective curve of genus $g\geq 3$ and let $X=C^g$ be the
$g$-copies of self products of $C$.
Then $\Griff_p(X)\otimes\Q$ are nontrivial for all $1\leq p\leq g-2$.
\end{proposition}

\medskip
The paper is organized as follows: \S2 is a review of the Lawson
homology groups, the duality between morphic cohomology and Lawson
homology, and intersection theory. In
 \S3 we discuss the $\mathbb{Q}$-coefficient Lawson homology and morphic
cohomology group of a finite quotient variety, and the intersection
theory in this setting. \S4 contains the main results that the
Lawson homology and morphic cohomology can be defined for the usual
Chow motives and the Chow motives of finite quotient projective
varieties. This is based on the fact that the action of
correspondences on Lawson homology are functorial, which is covered
in \S4.1 and \S4.2. As applications, in \S5.1 the projective bundle
theorem and blow-up formula for Lawson homology are reproved in a
different way; the computation of the Lawson homology for a smooth
projective variety with cell-decompositions is regained. \S5.2 is an
example on the finite quotient of projective variety: the
$\mathbb{Q}$-coefficient Lawson homology/morphic cohomology of
Hilbert schemes of points on a smooth surface. \S6 gives further
results and applications concerning generic rational maps, some new
examples with nontrivial Griffiths groups are built from known case.

\medskip\noindent\textbf{Acknowledgements.} We would like to thank
professor Blaine Lawson and Mark de Cataldo for suggestions,
conversations, encouragements and all their helps. The second author
also like to thank Jyh-Haur Teh for helpful conversations. We thank professor
E. Friedlander for his critical reading of the paper and numerous
suggestions and corrections.

\section{Lawson homology and morphic cohomology}
In this section, we first review the definition of the Lawson
homology $L_pH_k(X)$ for all integer $p$ and explain that the
properties of the original Lawson homology still
hold.\footnote{Friedlander pointed out to us that the consideration
for $p < 0$ is implicit in the work of Barry Mazur and himself, and
the formalism is worked out in the unpublished part of the thesis of
Mircea Voineagu.} Next, we will review morphic cohomology and the
duality.

\subsection{Lawson homology}

Denote by $H^{-1}\mathfrak{Abtop}$ the category of abelian
topological groups in which homotopy equivalences are inverted.

Given a projective variety $X$, we denote by ${\cal Z}_p(X)$ ($p\ge
0$) the space of algebraic $p$-cycles on $X$ with the natural Chow
topology. When $X$ is quasi-projective, Lima-Filho gave the
definition of ${\cal Z}_p(X)$ as the quotient
$${\cal Z}_p(X):={\cal Z}_p(\overline{X})/{\cal
Z}_p(\overline{X}-X),$$ where $\overline{X}$ is any projective
closure of $X$ (cf. \cite{Lima} and \cite{FG}). He shows that ${\cal
Z}_p(X)$ is well-defined up to isomorphism in the category
$H^{-1}\mathfrak{Abtop}$. As a consequence, the homotopy groups of
${\cal Z}_p(X)$ are independent of the choice of the projective
closure $\overline{X}$.

Based on the homotopy property of Lawson homology(\cite{FG},
Prop.2.3), the Lawson homology groups can be defined  for any
integer $p$ as follows, where $\mathbb{A}^r$ denotes the affine
 space of dimension $r$:
\begin{definition}\label{lawson homology}
Let $X$ be a quasi-projective variety. For a (possibly negative)
integer $p$, define the cycle space ${\cal Z}_p(X)$ to be the
homotopy equivalent class of ${\cal Z}_{p+r}(X\times \mathbb{A}^r)$
for any integer $r\ge \max(0,-p)$. (The homotopy property of Lawson
homology guarantees that ${\cal Z}_p(X)$ is independent of the
choice of $r$.)

The Lawson homology group $L_pH_k(X)$ is defined as
$$L_pH_k(X) :=\left\{
                \begin{array}{ll}
                  \pi_{k-2p}({\cal Z}_{p}(X)), & \hbox{if $ k\geq 2p$;} \\
                  0, & \hbox{if $k<2p$.}
                \end{array}
              \right.
 $$

\end{definition}

\begin{remark}
  For $p\ge 0$, the above definition coincides with the original definition of
Lawson homology groups. For $p<0$, we have
$L_pH_k(X)=\pi_{k-2p}({\cal Z}_{0}(X\times
\mathbb{A}^{-p}))=H^{BM}_{k-2p}(X\times
\mathbb{A}^{-p})=H^{BM}_k(X)=L_0H_k(X)$ (cf.
\cite{Friedlander-Haesemesyer-Walker}).
\end{remark}

Thus defined Lawson homology groups have expected functorial
properties.

\begin{definition}
  Let $f:X\to Y$ be a proper morphism between two quasi-projective
varieties. For $p\in\mathbb{Z}$, define the push-forward map $$f_*:
{\cal Z}_p(X)\to {\cal Z}_p(Y)$$ to be the one induced by
$$(f\times{\rm id})_*: {\cal Z}_{p+r}(X\times \mathbb{A}^r)\to {\cal
Z}_{p+r}(Y\times \mathbb{A}^r)$$ for a non-negative integer $r\ge
-p$.
\end{definition}

The following is essentially due to Frielander (\cite{Friedlander1}, Prop.2.9).
\begin{proposition}\label{functor push-forward}
  The above definition of $f_*$ does not depend on the choice of\;
$r$. Moreover, for two proper morphisms of quasi-projective
varieties $f: X\to Y$, $g: Y\to Z$, the following functoriality
holds:$$(gf)_*\stackrel{h.e.}{\simeq}g_*f_*: {\cal Z}_p(X)\to {\cal
Z}_p(Z), \quad \forall p\in\mathbb{Z}.$$
As a consequence,
$$(gf)_*=g_*f_*: L_pH_k(X)\to L_{p}H_{k}(Z), \quad \forall p, k\in\mathbb{Z}.$$
\end{proposition}

\medskip
Similarly, the definition of the pull-back map can be extended to
include the cycles of negative dimensions:
\begin{definition}\label{pull-back}
Let $X$ and $Y$ be quasi-projective varieties.  Let $f:X\to Y$ be a
l.c.i. (local complete intersection) morphism of codimension $d$
(i.e. $f$ factors into a regular imbedding followed by a smooth
morphism, and $\dim Y-\dim X=d$). For $p\in\mathbb{Z}$, define the
pull-back map in the category $H^{-1}\mathfrak{Abtop}$:
$$f^*: {\cal Z}_p(Y)\to {\cal Z}_{p-d}(X)$$ to be the one induced by
$$f^*: {\cal Z}_{p+r}(Y\times \mathbb{A}^r)\to {\cal
Z}_{p+r-d}(X\times \mathbb{A}^r)$$ for an integer $r\ge
\max(0,-p,d-p)$.
\end{definition}

\begin{proposition}\label{functor lawson}
  The above definition of $f^*$ does not depend on the choice of\;
$r$. Moreover, let $f: X\to Y$ and $g: Y\to Z$ be two l.c.i.
morphisms between quasi-projective varieties of codimension $d$ and
$e$ respectively. Then the following functoriality holds in the
category
$H^{-1}\mathfrak{Abtop}$:$$(gf)^*\stackrel{h.e.}{\simeq}f^*g^*:
{\cal Z}_p(Z)\to {\cal Z}_{p-d-e}(X), \quad \forall
p\in\mathbb{Z}.$$ As a consequence,
$$(gf)^*=f^*g^*: L_pH_k(Y)\to L_{p-d-e}H_{k-2d-2e}(X), \quad \forall p, k\in\mathbb{Z}.$$
\end{proposition}
\bp For the first part, it is enough to show the following diagram
commutes: $$\xymatrix{{\cal Z}_{p+r}(Y\times
\mathbb{A}^{r})\ar[r]^{(f\times{\rm id})^*}\ar[d]^{h.e.}_{\pi^*}&
{\cal
Z}_{p+r-d}(X\times \mathbb{A}^r)\ar[d]^{h.e.}_{\pi^*}\\
{\cal Z}_{p+r+1}(Y\times \mathbb{A}^{r+1})\ar[r]^{(f\times{\rm
id})^*}& {\cal Z}_{p+r+1-d}(X\times \mathbb{A}^{r+1})}$$ for $f$
being a smooth morphism and for $f$ being a regular imbedding,
respectively. In the former case we can check that the diagram
commutes by definition, the latter case is immediate from \cite{FG}
Theorem 3.4 (d).

To show the functoriality, choose an integer $r\ge \max(0, -p,
e-p,d+e-p)$ and consider
$$\xymatrix{{\cal Z}_{p+r}(Z\times \mathbb{A}^r)\ar[r]^{(f\times {\rm id})^*}& {\cal Z}_{p+r-e}(Y\times
\mathbb{A}^r) \ar[r]^{(g\times {\rm id})^*}& {\cal
Z}_{p+r-d-e}(X\times \mathbb{A}^r)}.$$ By the same method as in
\cite{Peters} Lemma 11c (also cf. proof of \cite{Fulton} Proposition
6.6(c)), we have $[(g\times {\rm id})\circ(f\times {\rm
id})]^*\stackrel{h.e.}{\simeq}(f\times {\rm id})^*(g\times {\rm
id})^*$. The conclusion follows.\qe

\begin{remark}
 If we define naively that $L_pH_k(X)=0$
for $p<0$, then the functoriality of pull-back maps does not hold in
general. For example, let $X=Z=\mathbb{P}^1$ and let $Y=pt$ be a
point. $f:X\to Y$ be the constant map and $g:Y\to Z$ maps $Y$ to any
point in $Z$. Then $$(g\circ f)^*: H_2(Z)=L_0H_2(Z)\to
L_0H_2(X)=H_2(X)$$ maps the generator of $H_2(Z)$ to the generator
of $H_2(X)$ by definition. On the other hand,
$$f^*g^*: L_0H_2(Z)\to L_{-1}H_0(Y)\to L_0H_2(X)$$ factor
through $L_{-1}H_0(Y)$. Thus, to save the functoriality ($(g\circ
f)^*=f^*g^*$), we have to define $L_pH_k(X)$ non-trivially.

We would like to point out here that the statement on ``$\alpha_*=0$
if $m<n-v$'' in Lemma 12 in Peters' paper \cite{Peters} is
imprecise. Despite of this minor imprecision,  his statement
$\alpha_*=0$ on $L_m^{hom}H_l(X)$ if $m<n-v$ is still valid. (cf.
Remark \ref{peters})
\end{remark}

\subsection{Morphic cohomology}

The morphic cohomology is defined by Friedlander and Lawson
\cite{FL1,FL2} by considering the homotopy groups of algebraic
cocycles.

Let $\mathfrak{Mor}(X,\mathcal{C}_r(Y))$ be the topological monoid
of effective algebraic cocycles of relative dimension $r$ with
values in $Y$, which by definition is the abelian monoid of
morphisms from $X$ to the Chow monoid $\mathcal{C}_r(Y)$ provided
with the compact open topology. When $X$ is geometrically unbranched
(e.g. when $X$ is normal), $\mathfrak{Mor}(X,\mathcal{C}_r(Y))$ can
be thought of as the subset of effective cycles in $X\times Y$ of
dimension $(r+\dim X)$ which is equidimensional over $X$.

The topological group $\mathcal{Z}^q(X)$ of all algebraic cocycles
of codimension $q$\, on $X$ is defined as a naive group completion
(cf. \cite{FL2} pg.538)
$$\mathcal{Z}^q(X):=[\mathfrak{Mor}(X,\mathcal{C}_0(\mathbb{P}^q))/
\mathfrak{Mor}(X,\mathcal{C}_0(\mathbb{P}^{q-1}))]^+=\mathfrak{Mor}(X,\mathcal{Z}_0(\mathbb{A}^q)).$$
Notice that the group $\mathcal{Z}^q(X)$ is not empty even when
$q\ge\dim X$.
\begin{definition} Let $X$ be a projective variety.\footnote{We restrict ourselves to consider the morphic cohomology of
 the projective varieties only, thanks to the reminder of Friedlander
that the formulation of morphic cohomology for quasi-projective
varieties is quite delicate.} The morphic cohomology groups are
defined to be the homotopy groups of $\mathcal{Z}^q(X)$:
$$L^qH^l(X):=\pi_{2q-l}(\mathcal{Z}^q(X)) \hbox{\quad if $2q\ge l$ and  $ q\ge 0$}.$$
If $2q<l$ or $q<0$, we define $L^qH^l(X)=0$.
\end{definition}

A duality map between morphic cohomology and Lawson homology is
observed by Friedlander and Lawson (\cite{FL2}). This duality can be
generalized with minor changes to include the case of any indices.
\begin{definition} 
  Let $X$ be a projective variety of dimension $m$. The duality map $$\mathcal{D}: \mathcal{Z}^p(X)\to
\mathcal{Z}_{m-p}(X), \quad \hbox{ $\forall p\le m$}£¬$$ is defined
by the graphing construction followed by the inverse of the natural
homotopy equivalence $\mathcal{Z}_{m-p}(X)\simeq
\mathcal{Z}_m(X\times \mathbb{A}^p)$ (which is meaningful even for
$p>m$ by Definition \ref{lawson homology}).

Taking the homotopy groups, we get the induced map (also denoted by
$\mathcal{D}$ by abuse of notation)
$$\mathcal{D}: L^pH^k(X)\to L_{m-p}H_{2m-k}(X).$$
\end{definition}

Now we recall some properties of the morphic cohomology groups and
the duality map that are needed in the rest of the paper.
\begin{enumerate}
  \item For any morphism $f: X\to Y$ between quasi-projective varieties, there is a pull-back morphism
   $$f^*: L^pH^k(Y)\to L^pH^k(X).$$
  When $f$ has equidimensional fibers (e.g., a flat morphism or a finite morphism) between normal varieties,
  there are Gysin homomorphisms $$f_!: L^pH^k(X)\to
L^{p-c}H^{k-2c}(Y)$$ for $2p\ge k\ge 2c$, where $c=\dim X-\dim Y$
(implied by \cite{FL1} Proposition 2.5).

  When $X$ and $Y$ are smooth projective varieties,
let $c=\dim X-\dim Y$ and define $$f_*: L^pH^k(X)\to
L^{p-c}H^{k-2c}(Y),\quad \forall p,k\in\mathbb{Z}$$ as
$\mathcal{D}^{-1} f_*\mathcal{D}$, where $\mathcal{D}$ is the
duality between morphic cohomology and Lawson homology defined in
Proposition \ref{negative} and $f_*$ is the push-forward for Lawson
homology.

Similarly, given $f': X'\to Y'$ between finite quotient of smooth
projective varieties, we define
$$f'_*:=(\mathcal{D}')^{-1} f'_*\mathcal{D}':
L^pH^k(X',\mathbb{Q})\to L^{p-c}H^{k-2c}(Y',\mathbb{Q})$$ where let
$c=\dim X'-\dim Y'$, the push-forward map $f'_*$ on the right hand
side is the one for Lawson cohomology and
$\mathcal{D}'=\mathcal{D}\otimes\mathbb{Q}$ is induced by the
duality map $\mathcal{D}$ (which is an isomorphism by Lemma \ref{2
diagrams}).

 \item There is a cup product $$\#: L^pH^k(X)\otimes L^{p'}H^{k'}(X)\to L^{p+p'}H^{k+k'}(X)$$
  natural with respect to morphisms (\cite{FL1} Corollary 6.2), i.e.
for $f:X\to Y$ a morphism between quasi-projective varieties,
  $$f^*(\alpha\#\beta)=f^*(\alpha)\#f^*(\beta), \quad \forall \alpha, \beta\in L^*H^*(Y).$$

  \item
If $X$ is smooth and projective, then the duality map
$$\mathcal{D}: L^pH^k(X)\to L_{m-p}H_{2m-k}(X), \hbox{ for $p\le m$ }$$ is an isomorphism
  compatible with the ring structures of $L^*H^*(X)$ and $L_*H_*(X)$,
  i.e.
  $$\mathcal{D}(\alpha\#\beta)=D(\alpha)\bullet D(\beta), \quad\forall \alpha\in L^pH^k(X),\beta
\in L^qH^l(X),$$ where $p,q\le m$ and $p+q\le m$. (This restriction
on $p, q$ is unnecessary, see Proposition \ref{negative} below.)
\end{enumerate}

The duality behaves as expected for the Lawson homology with
possibly negative dimension:
\begin{proposition}\label{negative} (1) Let $f: X\to Y$ be a morphism between projective
varieties with dimension $m$ and $n$, respectively. Then
$f^*\mathcal{D}=\mathcal{D}f^*$. In another word, the following
diagram commutes for any integer $p$:
$$\xymatrix{L^pH^k(Y)\ar[r]^-{\mathcal{D}}\ar[d]^{f^*}& L_{n-p}H_{2n-k}(Y)\ar[d]^{f^*}\\ L^pH^k(X)\ar[r]^-{\mathcal{D}}& L_{m-p}H_{2m-k}(X)} $$

(2) If $X$ is smooth and projective, then the duality (which we call the  Friedlander-Lawson duality)
$$\mathcal{D}: L^pH^k(X)\to L_{m-p}H_{2m-k}(X)$$ is a group isomorphism for any integer $p$.

(3) $\mathcal{D}(\alpha\#\beta)=D(\alpha)\bullet D(\beta), \forall
\alpha,\beta
  \in L^*H^*(X).$ In another word, the duality map is compatible with the ring structures of morphic cohomology and of Lawson
homology.
\end{proposition}
\bp (1) is an immediate consequence of the known result (cf.
\cite{FL2} Proposition 2.2 and 2.3).

(2) follows directly from the proof of \cite{FL2} Theorem 3.3.

(3) is proved in the next section \S\ref{intersection}.\qe

\subsection{Intersection theory}\label{intersection}

 In this section, assume $X$ is a smooth quasi-projective
variety. Then the diagonal map $\Delta: X\to X\times X$ is a regular
imbedding. Let $\alpha\in L_pH_k(X), \beta\in L_qH_l(X)$. There is a
natural map $\mathcal{Z}_p(X)\wedge \mathcal{Z}_q(X)\to
\mathcal{Z}_{p+q}(X\times X)$, where ``$\wedge$'' is the smash
product. Taking the homotopy groups at both sides, we get a natural
map $L_pH_k(X)\times L_qH_l(X)\to L_{p+q}H_{k+l}(X\times X)$, and we
denote the image of $(\alpha,\beta)$ under this map by
$\alpha\times\beta$.

\begin{definition}\label{bullet}
  Let $X$ be a smooth quasi-projective variety and
   $\Delta: X\to X\times X$ be the diagonal map.
For any $\alpha\in L_pH_k(X), \beta\in L_qH_l(X)$, the intersection
$\alpha\bullet \beta\in L_{p+q-m}H_{k+l-2m}(X)$ is defined as
$$\alpha\bullet\beta:=\Delta^*(\alpha\times \beta).$$
\end{definition}

Notice that in the above definition, no restriction is put on
$p,q,k,l$ and fortunately, the compatibility with pull-back
$f^*(\alpha\bullet\beta)=f^*\alpha\bullet f^*\beta$, the
compatibility with duality
$\mathcal{D}(\alpha\#\beta)=D(\alpha)\bullet D(\beta)$ and the
projection formula $f_*(f^*\alpha\bullet\beta)=\alpha\bullet
f_*\beta$ still hold in this more general situation where the cycles
of negative dimensions are allowed. The proof are essentially the
same as the canonical case. We explain as follows.

\medskip

 First we prove the compatibility with duality:

\bp(of Proposition \ref{negative} (3):
$\mathcal{D}(\alpha\#\beta)=D(\alpha)\bullet D(\beta)$.) The proof
is exactly the same as \cite{FL2} Proposition 2.7 and its remark,
where cycle spaces of negative dimensions, if appear, are understood
as in Definition \ref{lawson homology}.
 \qe

\medskip
The next proposition (1)(2)(3) is adapted from \cite{Peters} Lemma
11, with minor revises, while (4) is adapted from \cite{FG} Theorem
3.5 (b).
\begin{proposition}\label{3 formulas} Let $X, Y, X', Y'$ be smooth quasi-projective
varieties.

(1) Let $f: X\to Y$ be a morphism. Then
$$f^*(\alpha\bullet\beta)=f^*\alpha\bullet f^*\beta, \quad \forall
\alpha\in L_pH_k(Y), \beta\in L_qH_l(Y).$$

(2)Suppose the following is a fibre square: $$\xymatrix{X'\ar[r]^{f'}\ar[d]^{g'}&Y'\ar[d]^{g}\\
X\ar[r]^f& Y}$$ Then with $d=\dim Y-\dim X$, for any integers $p,k$,
one has
$$f^*g_*=g'_*f'^*: L_pH_k(Y')\to L_{p-d}H_{k-2d}(X).$$

(3) Let $f:X\to Y$ be a morphism. Then the projection formula holds
in for any integers $p,q,k,l$:
$$f_*(\alpha\bullet f^*\beta)=f_*\alpha\bullet\beta, \quad \alpha\in L_pH_k(X), \beta\in L_qH_l(Y).$$

(4) The intersection is graded-commutative and associative, i.e.,
$\forall\alpha\in L_pH_k(X),\beta\in L_qH_l(X),\gamma\in L_rH_m(X)$,
we have
$$\alpha\bullet\beta=(-1)^{kl}\beta\bullet\alpha,$$
$$(\alpha\bullet\beta)\bullet\gamma=\alpha\bullet(\beta\bullet\gamma).$$
\end{proposition}
\bp (1) It is a immediate consequence of functoriality of pull-back
(Proposition \ref{functor lawson}). Indeed, since any morphism
between smooth quasi-projective varieties are l.c.i., we have
$$f^*(\alpha\bullet\beta)=f^*\Delta_X^*(\alpha\times\beta)=\Delta_Y^*(f\times f)^*(\alpha\times\beta)=f^*\alpha\bullet f^*\beta.$$

(2) Take an integer $r\ge \max(0, -p, -p+d)$ and consider the
following diagram in the category $H^{-1}\mathfrak{Abtop}$:
$$\xymatrix{\mathcal{Z}_{p-d+r}(X'\times \mathbb{A}^r)\ar[d]^{(g'\times{\rm id})_*}&\mathcal{Z}_{p+r}(Y'\times \mathbb{A}^r)\ar[d]^{(g\times{\rm id})_*}\ar[l]_-{(f'\times{\rm id})^*}\\
\mathcal{Z}_{p-d+r}(X\times \mathbb{A}^r)& \mathcal{Z}_{p+r}(Y\times
\mathbb{A}^r)\ar[l]_-{(f\times{\rm id})^*}}$$ It commutes, by
consider the case when $f$ is a regular imbedding (\cite{FG} Theorem
3.4 d) and the case when $f$ is a flat morphism (\cite{Fulton}
Proposition 1.7). Then by our definition of cycle spaces (Definition
\ref{lawson homology}), the conclusion follows.

(3)  Peters' proof is valid in this setting: let $\gamma_f=({\rm
id}, f)=({\rm id}\times f)\circ \Delta_X: X\to X\times Y$. Consider
the fibre square
$$\xymatrix{X\ar[r]^-{\gamma_f}\ar[d]^{f}& X\times Y\ar[d]^{f\times\textrm{id}}\\
Y\ar[r]^-{\Delta_Y}&Y\times Y}$$ Then $f_*\alpha\bullet
\beta=\Delta^*_Y(f_*\beta\times \beta)=\Delta^*_Y(f\times
\textrm{id})_*(\alpha\times \beta)=f_*((\gamma_f)^*(\alpha\times
\beta))=f_*(\Delta_X^*({\rm id}\times f)^*(\alpha\times
\beta))=f_*(\Delta_X^*(\alpha\times f^*\beta))=f_*(\alpha\bullet
f^*\beta)$, where the third equality is by applying (2) to the above
fibre square, the fourth is by the functoriality of pull-back
(Proposition \ref{functor lawson}).

(4) The standard argument still applies for negative cycle
spaces.\qe

\section{Quotient variety by a finite group}

In this section we explore the relation between the Lawson homology
of a quasi-projective variety and the Lawson homology of its finite
quotient. The goal is to establish Proposition \ref{G-invariant}.

Suppose a finite group $G$  acts faithfully on a quasi-projective
variety $X$ (that is, the only element in $G$ fixing every point in
$X$ is the identity). The quotient $X'=X/G$ is again a
quasi-projective variety (cf. \cite{Harris} \S10). Let $\pi:X\to X'$
denote the quotient map. We give the definition of pull-back $\pi^*$
of algebraic cycles as in \cite{Fulton} Example 1.7.6:

For any subvariety $W$ of $X$, let
$$I_W=\{g\in G:    g|_W=\hbox{id}_W\}$$ be the inertia group. Let
$e_W=\hbox{card}(I_W)$ be the order of the group $I_W$.

\begin{definition}\label{pi*} A group homomorphism $\pi^*: \mathcal{Z}_p(X')\to \mathcal{Z}_p(X)$
for $p\ge 0$ is defined as follows: for any subvariety $V$ of $X'$,
let
$$\pi^*[V]=\sum e_W[W],$$
the sum over all irreducible components of $\pi^{-1}(V)$.

In general, for a possibly negative $p$ we take $r\ge \max(0,-p)$
and define $\pi^*:
\mathcal{Z}_{p+r}(X'\times\mathbb{A}^r)\to\mathcal{Z}_{p+r}(X\times\mathbb{A}^r)$
same as above, which induces $\pi^*: \mathcal{Z}_p(X')\to
\mathcal{Z}_p(X)$.
\end{definition}

\begin{remark} This definition is uniquely characterized
by the fact that $$\pi^*\pi_*[W]=G[W]:=\sum_{g\in G} g_*[W].$$
\end{remark}

To induce from $\pi^*$ a map between Lawson homology groups of $X'$
and $X$, it is necessary to verify the continuity of $\pi^*$.

\begin{lemma}\label{continuity}
 Let $\pi:X\to X':=X/G$ denote the quotient
 map of a quasi-projective variety with a faithful action of a finite group $G$. The map $\pi^*: \mathcal{Z}_p(X')\to \mathcal{Z}_p(X)$ is
 continuous. As a consequence, it induces a morphism $\pi^*: L_pH_k(X')\to
 L_pH_k(X)$.
\end{lemma}
\bp Without loss of generality, we can assume $p\ge 0$, since the
case when $p<0$ can be easily deduced from the case $p=0$.

  By \cite{Lima2} Theorem 3.1 and Theorem 5.8, for a complex algebraic variety (in particular, a
complex quasi-projective variety) $X$, there are three equivalent
definitions for the topology of $\mathcal{Z}_p(X)$, namely, the
  flat topology $\mathcal{Z}_p(X)^{fl}$, the equidimensional topology $\mathcal{Z}_p(X)^{eq}$, and
  $\mathcal{Z}_p(X)^{ch}$ defined
  via Chow varieties (which is the original definition of the topology of $\mathcal{Z}_p(X)$). Therefore it suffices to show the
  continuity for
  $$\pi^*: \mathcal{Z}_p(X')^{fl}\to \mathcal{Z}_p(X)^{eq}.$$
  Let $S$ be a smooth projective variety. Given a cycle $\Gamma'$ on $S\times
  X'$ which is flat over $S$ and of relative dimension $p$. Consider
  $\Gamma:=(\hbox{id}\times \pi)^*(\Gamma')$ (which is well defined since
  $\hbox{id}\times \pi$ is a finite quotient morphism). Notice that
  $\Gamma$ may not be flat over $S$, but is still equidimensional over $S$
  of relative dimension $p$. Given $s\in S$, let
  $[\Gamma_s]$ be the intersection theoretic fiber over $s$. Then it
  suffices to show that $\pi^*([\Gamma'_s])=[\Gamma_s]$ for any $s\in S$.

Notice that  $(\hbox{id}\times \pi)_*[\Gamma]=(\hbox{id}\times
\pi)_*(\hbox{id}\times \pi)^*[\Gamma']=|G| [\Gamma']$. Then
\begin{eqnarray*}
\pi_*[\Gamma_s]&=&(\hbox{id}\times \pi)_*([\Gamma]\cdot [s\times X])\\
  &=&(\hbox{id}\times\pi)_*\big{(}[\Gamma]\cdot(\hbox{id}\times\pi)^*[s\times
  X']\big{)}\\
  &=&(\hbox{id}\times\pi)_*[\Gamma]\cdot[s\times X']\\
  &=&|G| [\Gamma'_s].
\end{eqnarray*}
The notation $\cdot$ denotes the refined intersection. The third
equality is because of the projection formula for refined
intersection.

Therefore
$$|G|\pi^*([\Gamma'_s])=\pi^*(|G|[\Gamma'_s])=\pi^*\pi_*[\Gamma_s]=G[\Gamma_s]=|G|\,[\Gamma_s],$$
where the last equality is by the invariance of $\Gamma_s$ under the
action of $G$.

Since $\pi^*$ is a morphism of free abelian groups, so by dividing
$|G|$ from both sides of the above equalities we conclude that
$$\pi^*(\Gamma_s')=\Gamma_s,$$ which completes the proof.
 \qe

\medskip

We need the following elementary fact about homotopy groups of
topological abelian groups.

\begin{lemma}\label{lemma homotopy}
{\label{sec:lemma2.1} Let $f_1, f_2: Z_1\rightarrow Z_2$ be two
continuous homomorphisms between topological abelian groups. Then
the induced homomorphism of the sum is the sum of the induced
homomorphisms on the homotopy groups, i.e.,
$(f_1+f_2)_*=(f_1)_*+(f_2)_*:\pi_k(Z_1)\rightarrow \pi_k(Z_2)$, for
any integer $k\ge 0$. }
\end{lemma}

\bp Let $\alpha\in \pi_k(Z_1)$ and $g\in \alpha$. Sometimes we also
write  $[g]=\alpha$. Then
$(f_1+f_2)_*(\alpha)=(f_1+f_2)_*([g])=[(f_1+f_2)\circ g]=[f_1\circ
g+f_2\circ g]=[f_1\circ g]+[f_2\circ
g]=(f_1)_*([g])+(f_2)_*([g])=(f_1)_*(\alpha)+(f_2)_*(\alpha)$. That
is what we want to prove. \qe

\medskip Now we show the following relation between the Lawson
homology groups of $X$ and its quotient $X'=X/G$.
\begin{proposition}
\label{G-invariant}
Let $\pi:X\to X':=X/G$ denote the quotient
 map of a quasi-projective variety  with a faithful action of a finite group $G$. Then there is a
canonical isomorphism
$$\pi_*: (L_pH_k(X,\mathbb{Q}))^G\cong L_pH_k(X',\mathbb{Q}), \quad \hbox{ for any $p, k\in \mathbb{Z}$}.$$
 and  an isomorphism (if $X$ is projective)
$$\pi_!: (L^qH^l(X,\mathbb{Q}))^G\cong L^qH^l(X',\mathbb{Q}), \quad \hbox{ for any $q, l\in \mathbb{Z}$}.$$
\end{proposition}

\bp We provide here the proof of the isomorphism of $\pi_*$, since
the isomorphism of $\pi_!$ can be proved similarly.

Consider the push-forward map $\pi_*$ and the pull-back $\pi^*$
which is continuous by Lemma \ref{continuity}. It is easy to verify
from the definition that, on the cycle spaces,
 $$\pi_*\pi^*=|G| \cdot \hbox{id}:{\mathcal{Z}}_p(X')\rightarrow
  {\mathcal{Z}}_p(X')
  $$and
  $$\pi^*\pi_*=\sum_{g\in G}g_*: {\mathcal{Z}}_p(X)\rightarrow
  {\mathcal{Z}}_p(X) .
   $$
Therefore, we have corresponding identities on Lawson homology
groups, by a property of homotopy groups of topological abelian
groups (Lemma \ref{lemma homotopy}):

$$\pi_*\pi^*=|G| \cdot \hbox{id}: L_pH_k(X')\to
L_pH_k(X'),
$$ and
$$\pi^*\pi_*=\sum_{g\in G}g_*: L_pH_k(X)\to
L_pH_k(X).
   $$ Then the conclusion follows by the following simple fact about vector
   spaces (Lemma \ref{lemma isom}).
\qe

   \begin{lemma}\label{lemma isom}
Let $V_1$, $V_2$ be two $\mathbb{Q}$-vector spaces acted by a finite
group $G$. Suppose $G$ acts trivially on $V_2$ and denote the $G$-invariant subspace
of $V_1$ by $V_1^G$. Let $\phi: V_1\to V_2$,
$\psi: V_2\to V_1$ be two equivariant linear maps of vector spaces
(i.e. $\phi(gx)=\phi(x)$ and $g\psi(y)=\psi(y)$, $\forall g\in G, x\in V_1, y\in V_2$).
 If the following two conditions are satisfied,
 \begin{itemize}
 \item[i)]  $\forall y\in V_2$, $\phi\circ \psi(y)=|G|\cdot y$,
 \item[ii)] $\forall x\in V_1$, $\psi\circ \phi(x)=Gx:=\sum_{g\in
 G}gx$,\end{itemize}
  then $\phi|_{V_1^G}: V_1^G\to V_2$ is an isomorphism, with inverse
  $\psi/ |G|$.
\end{lemma}
\bp The surjectivity of $\phi|_{V_1^G}$ is because of the
surjectivity of $\phi\circ \psi=|G|\cdot \hbox{id}_{V_2}$. For
injectivity, suppose $x\in V_1^G$ satisfying $\phi(x)=0$ Since $x$
is invariant under $G$-action, $0=\psi\circ \phi(x)=Gx=|G|\cdot x$,
which implies that $x=0$.
 \qe

\medskip

Next, we define a natural intersection ring structure on the
$\mathbb{Q}$-coefficient Lawson homology groups of a finite quotient
of a smooth quasi-projective variety.

\begin{definition}\label{dot} Let $X$ be a smooth quasi-projective variety with a finite
group $G$ acting on it faithfully. Denote the quotient map by $\pi:
X\to X'$. For any $\alpha\in L_pH_k(X',\mathbb{Q})$ and $\beta\in
L_qH_l(X',\mathbb{Q})$, the intersection $\alpha\cdot\beta$ in
$L_{p+q-m}H_{k+l-2m}(X',\mathbb{Q})$ is defined as
\begin{equation}
\alpha\cdot\beta:=\frac{1}{|G|}\pi_*(\pi^*\alpha\bullet\pi^*\beta)
\end{equation}
where $\pi^*$ is defined in Definition \ref{pi*}, and $\bullet$ is
defined in Definition \ref{bullet}.
\end{definition}

\begin{proposition}\label{independent}
  Assume further that $X'$, hence $X$, is projective. Then the intersection product defined as above
  depends only on $X'$, not on the choice of $X$
  and  $G$.
\end{proposition}
The proof is postponed to the end of this section. Our method is to
 compare the above intersection product with the cup product of the morphic cohomology.

\begin{lemma} Use the notation as in the above Definition \ref{dot}. Then for
any $p,q,r,k,l,m\in\mathbb{Z}$, $\alpha\in L_pH_k(X',\mathbb{Q})$,
$\beta\in L_qH_l(X',\mathbb{Q})$  and $\gamma\in
L_rH_m(X,\mathbb{Q})$, we have

(1) $\pi^*(\alpha\cdot\beta)=\pi^*(\alpha)\bullet\pi^*(\beta)$,

(2) $\pi_*((\pi^*\alpha)\bullet \gamma)=\alpha\cdot\pi_*(\gamma)$.
\end{lemma}
\bp (1) The definition of $\pi^*$ and Proposition \ref{3 formulas}
imply that both sides are invariant under the $G$-action. Therefore
it is enough to show the equality
$$\pi_*\pi^*(\alpha\cdot\beta)=\pi_*(\pi^*(\alpha)\bullet\pi^*(\beta)).$$
Since $\pi_*\pi^*=|G|\cdot{\rm id}$, then by Definition \ref{dot}
the above equality holds.

(2) The right hand side equals to
$$\frac{1}{|G|}\pi_*(\pi^*\alpha\bullet \pi^*\pi_*\gamma)=
\frac{1}{|G|}\pi_*(\pi^*\alpha\bullet \sum_{g\in G} g^*\gamma)
=\frac{1}{|G|}\sum_{g\in G}\pi_*g^*((g^{-1})^*\pi^*\alpha\bullet
\gamma).$$ Since $\pi^*\alpha$ is $G$-invariant,
$(g^{-1})^*\pi^*\alpha=\pi^*\alpha$. Moreover, $\pi_*g^*=\pi_*$.
Therefore the above equals to the left hand side
$\pi_*((\pi^*\alpha)\bullet \gamma)$. \qe

\begin{lemma}\label{2 diagrams} Let $\pi: X\to X'=X/G$ be a finite quotient map where $G$ acts
faithfully on a projective normal variety $X$, then the following
diagrams commute (here we denote by $\mathcal{D}'$ the duality map
for $X'$):
\begin{equation}\label{diagram1}
\xymatrix{L^rH^k(X,\mathbb{Q})\ar[r]^-{\mathcal{D}}\ar[d]^{\pi_!} &
L_{m-r}H_{2m-k}(X,\mathbb{Q})\ar[d]^{\pi_*} \\
L^rH^k(X',\mathbb{Q})\ar[r]^-{\mathcal{D}'}&L_{m-r}H_{2m-k}(X',\mathbb{Q})}
\end{equation}

$$\xymatrix{
L^rH^k(X,\mathbb{Q})\ar[r]^-{\mathcal{D}}&L_{m-r}H_{2m-k}(X,\mathbb{Q}) \\
L^rH^k(X',\mathbb{Q})\ar[r]^-{\mathcal{D}'}\ar[u]^{\pi^*}&
L_{m-r}H_{2m-k}(X',\mathbb{Q})\ar[u]^{\pi^*}}$$ As a consequence,
$\mathcal{D}'$ is an isomorphism. \end{lemma} \bp It is easy to
check that both diagrams hold on the level of cocycles (in place of
morphic cohomology) and cycles (in place of Lawson homology). Then
by Proposition \ref{G-invariant} and the fact that $\mathcal{D}$ is
an isomorphism, we know $\mathcal{D}'$ is also an isomorphism.\qe

\begin{proposition}\label{duality}
  Let $\pi: X\to X'=X/G$ be a finite quotient map where $G$ acts
faithfully on a smooth projective variety $X$. Then for any
$\alpha,\beta\in
  L^*H^*(X',\mathbb{Q})$, the duality map $\mathcal{D}': L^*H^*(X',\mathbb{Q})\to
  L_*H_*(X',\mathbb{Q})$ satisfies
  $$\mathcal{D}'(\alpha\#\beta)=\mathcal{D}'(\alpha)\cdot \mathcal{D}'(\beta).$$
\end{proposition}
\bp We have
$$\begin{array}{lll}
|G|\mathcal{D}'(\alpha)\cdot\mathcal{D}'(\beta)&=
\pi_*\big{[}\pi^*\mathcal{D}'(\alpha)\bullet\pi^*\mathcal{D}'(\beta)\big{]}
=\pi_*\big{[}\mathcal{D}(\pi^*\alpha)\bullet\mathcal{D}(\pi^*\beta)\big{]}\\
&=\pi_*\mathcal{D}(\pi^*\alpha\#\pi^*\beta)=\pi_*\mathcal{D}
\pi^*(\alpha\#\beta)=\mathcal{D}'\pi_*\pi^*(\alpha\#\beta)=|G|\mathcal{D}'(\alpha\#\beta).
\end{array}$$ where the second and fifth equalities are because of
Lemma \ref{2 diagrams}, the third is from Proposition \ref{negative}
(3), the fourth is because the pull-back $\pi^*$ is compatible with
the product of morphic cohomology. \qe

\medskip
\bp (of Proposition \ref{independent}) By the above proposition, it
is enough to show that
 the duality map $\mathcal{D}'$ is surjective, since then the product in the Lawson homology $L_*H_*(X',\mathbb{Q})$
  is determined by the cup product $\#$ in the morphic cohomology $L^*H^*(X',\mathbb{Q})$.

On the other hand, by assumption $X$ is smooth projective, then
$\mathcal{D}$ is an isomorphism. Lemma \ref{sec:lemma2.1}
asserts that $\pi_*$ is an isomorphism. Then by diagram
(\ref{diagram1}) we know $\mathcal{D}'\pi_!=\pi_*\mathcal{D}$ is
surjective, it follows that $\mathcal{D}'$ is surjective.
  \qe



\section{Correspondences and Motives}\label{sec:2}

\subsection{The action of Correspondences between smooth varieties} Let $X$ and $Y$ be smooth
projective varieties. A \textbf{correspondence} $\Gamma$ from $X$ to
$Y$ is a cycle (or an equivalent class of cycles depending on the
context) on $X\times Y$.
%
%
We denote the group of correspondences of rational equivalence
classes between varieties $X$ and $Y$ by
$$Corr_d(X, Y):={\rm Ch}_{\dim X+d}(X\times Y).$$
In general without assuming the varieties $X, Y$ to be connected, we
define $$Corr_d(X, Y):=\oplus{\rm Ch}_{\dim
X_\alpha+d}(X_\alpha\times Y),$$  where $X=\coprod X_\alpha$ is the
decomposition of connected components of $X$.

\medskip

Recall (\cite{Fulton}, Chapter 16) that a correspondence $\Gamma\in
Corr_d(X, Y)$ acts on Chow groups as follows
$$
\begin{array}{cc}
&\Gamma_*: {\rm Ch}_p(X)\rightarrow {\rm Ch}_{p+d}(Y)\\
& \Gamma_*(u)=p_{2*}(p_1^*u\bullet \Gamma)
\end{array}
$$
where $p_1$ (resp. $p_2$) denote the projection from $X\times Y$
onto $X$ (resp. $Y$) and  $\bullet$ is the intersection product on
the Chow group of the smooth variety $X\times Y$.

Let $X$, $Y$, $Z$ be smooth projective varieties. The composition of
two correspondences $\Gamma_1\in Corr_{d_1}(X,Y)$ and $\Gamma_2\in
Corr_{d_2}(Y,Z)$ is given by the formula
$$\Gamma_2\circ\Gamma_1=p_{13*}(p_{12}^*\Gamma_1\cdot p_{23}^*\Gamma_2)\in Corr_{d_1+d_2}(X, Z)
$$
where $p_{ij}$, $i,j =1, 2, 3$ are the projection of $X\times Y\times Z$ on the
product of its $i$th and $j$th factors.


Follow the idea of Peters \cite{Peters}, we define the analogous
homomorphisms on the level of Lawson homology by the same formula.
Notice that for any $\Gamma\in Corr_d(X, Y)$,  by modulo algebraic
equivalence instead of rational equivalence relation it determines
an element in
$$\pi_0(\mathcal{Z}_{\dim X+d}(X\times Y))=L_{\dim X+d}H_{2\dim
X+2d}(X\times Y)$$ which is again denoted by $\Gamma$ by abuse of
notation.

\begin{definition}\label{push-forward}
Let $X$, $Y$ be smooth projective varieties, $\Gamma\in Corr_d(X,
Y)$. Then for any element $\alpha\in L_pH_k(X)$, the push-forward
morphism is defined by
$$
\begin{array}{cc}
&\Gamma_*: L_pH_k(X)\rightarrow L_{p+d}H_{k+2d}(Y)\\
& \Gamma_*(\alpha)=p_{2*}(p_1^*\alpha\bullet \Gamma).
\end{array}
$$
\end{definition}

\medskip

Now we show that the push-forward morphism defined as above is
functorial.

\begin{proposition}\label{composition}
  Let $X, Y, Z$ be smooth projective varieties, $\Gamma_1\in
  Corr_d(X,Y)$ and $\Gamma_2\in Corr_e(Y, Z)$. Then for any $u\in
  L_pH_k(X)$, we have
  $$(\Gamma_2\circ\Gamma_1)_*u=\Gamma_{2*}\Gamma_{1*}u \in L_{p+d+e}H_{k+2d+2e}(Z).$$
\end{proposition}
\bp The proof is by applying basic properties of push-forward and
pull-back of Lawson homology groups which we list below for
convenience:

\begin{enumerate}
  \item Graded commutativity and associativity (Proposition \ref{3 formulas} (4)).
  \item Functoriality of push-forward and pull-back: $(fg)^*=g^*f^*$, $(fg)_*=f_*g_*$ (Proposition \ref{functor push-forward} and \ref{functor lawson}).
  \item Projection formula: $f_*(\alpha\bullet
f^*\beta)=f_*\alpha\bullet\beta$ (Proposition \ref{3 formulas} (3)).
  \item Pull-back compatible with the intersection product: $f^*(\alpha\bullet\beta)=f^*\alpha\bullet f^*\beta$  (Proposition \ref{3 formulas} (1)).
  \item Given a fiber square $$\xymatrix{W\ar[r]^g\ar[d]^q&Z\ar[d]^p\\Y\ar[r]^f&X}$$ where $f, g$ are proper, and  $p, q$ are
  flat,
  then $p^*f_*=g_*q^*$ (Proposition \ref{3 formulas} (2)).

\end{enumerate}

Denote by $p^{XYZ}_{XY}$ the projection from $X\times Y\times Z$ to
$X\times Y$, and similarly for other projections.

$$\begin{array}{lll}(\Gamma_2\circ\Gamma_1)_*u&=p^{XZ}_{Z*}\big{(}(\Gamma_2\circ\Gamma_1)\bullet p^{XZ*}_Xu\big{)}\\
&=p^{XZ}_{Z*}\big{(}p^{XYZ}_{XZ*}(p^{XYZ*}_{XY}\Gamma_1\bullet p^{XYZ*}_{YZ}\Gamma_2)\bullet p^{XZ*}_Xu\big{)}\\
&\stackrel{3}{=}p^{XZ}_{Z*}p^{XYZ}_{XZ*}\big{(}(p^{XYZ*}_{XY}\Gamma_1\bullet p^{XYZ*}_{YZ}\Gamma_2)\bullet p^{XYZ*}_{XZ}p^{XZ*}_Xu\big{)}\\
&\stackrel{2}{=}p^{XYZ}_{Z*}\big{(}(p^{XYZ*}_{XY}\Gamma_1\bullet p^{XYZ*}_{YZ}\Gamma_2)\bullet p^{XYZ*}_{X}u\big{)}\\
&\stackrel{1}{=}p^{XYZ}_{Z*}\big{(}p^{XYZ*}_{YZ}\Gamma_2\bullet(p^{XYZ*}_{XY}\Gamma_1\bullet p^{XYZ*}_{X}u)\big{)}\\
&\stackrel{2,4}{=}p^{YZ}_{Z*}p^{XYZ}_{YZ*}\big{(}p^{XYZ*}_{YZ}\Gamma_2\bullet(p^{XYZ*}_{XY}(\Gamma_1\bullet p^{XY*}_{X}u)\big{)}\\
&\stackrel{3}{=}p^{YZ}_{Z*}\big{(}\Gamma_2\bullet p^{XYZ}_{YZ*}p^{XYZ*}_{XY}(\Gamma_1\bullet p^{XY*}_{X}u)\big{)}\\
&\stackrel{5}{=}p^{YZ}_{Z*}\big{(}\Gamma_2\bullet p^{YZ*}_{Y}p^{XY}_{Y*}(\Gamma_1\bullet p^{XY*}_{X}u)\big{)}\\
&=\Gamma_{2*}\Gamma_{1*}u.
\end{array}
$$ where the first and last equalities hold by the definition of push-forward for Lawson homology (Definition
\ref{push-forward}). For the second equality we use the \cite{FG}
Theorem 3.5 c, which asserts that for an intersection pairing
$\mathcal{Z}_p(X)\times \mathcal{Z}_q(X)\stackrel{\bullet}{\to}
\mathcal{Z}_{p+q-\dim X}(X)$, applying $0$-th homotopy $\pi_0$
yields the usual intersection product on algebraic equivalence
class, hence is compatible with the ring structure for Chow groups.
\qe

\medskip

Denote $Corr^{d}(X,Y):=Ch^{\dim X+d}(X\times Y)=Ch_{\dim
Y-d}(X\times Y)$. By the duality isomorphism $\mathcal{D}$ between
morphic cohomology and Lawson homology, the analogous functorial
property for morphic cohomology immediately follows:

\begin{proposition}\label{composition2}
  Let $X, Y, Z$ be smooth projective varieties, $\Gamma_1\in
  Corr^{d}(X,Y)$ and $\Gamma_2\in Corr^{e}(Y, Z)$. Then for any $u\in L^qH^l(X)$,
  $$(\Gamma_2\circ\Gamma_1)_*u=\Gamma_{2*}\Gamma_{1*}u \in L^{q+d+e}H^{l+2d+2e}(Z).$$
\end{proposition}
\qe

\subsection{The action of Correspondences between quotient varieties}
In this subsection, we extend the action of correspondences to the
category of finite quotients of nonsingular projective varieties.
The definition is formally the same as Definition
\ref{push-forward}, with the intersection as defined in Definition
\ref{dot}.

\begin{definition}
Let $X'$, $Y'$ be two finite quotient varieties and let $\Gamma'\in
Corr_d(X', Y')_{\mathbb{Q}}$. The push-forward $\Gamma'_*$ is
defined by
$$
\begin{array}{cc}
&\Gamma'_*: L_pH_k(X',\mathbb{Q})\rightarrow L_{p+d}H_{k+2d}(Y',\mathbb{Q})\\
& \Gamma'_*(u)=p_{2*}(p_1^*u\cdot \Gamma')
\end{array}
$$\end{definition}
Note that the above definition implicitly uses Proposition
\ref{independent}, i.e. the intersection product on a finite
quotient variety is well-defined.

An important property for the push-forward action of a
correspondence is the following functoriality.
\begin{proposition}\label{quotient functorial}
Let $X', Y', Z'$ be finite quotient varieties, $\Gamma'_1\in
  Corr_d(X',Y')$ and $\Gamma'_2\in Corr_e(Y', Z')$. Then for any $u\in L_pH_k(X',\mathbb{Q})$,
  $$(\Gamma'_2\circ\Gamma'_1)_*u=\Gamma'_{2*}\Gamma'_{1*}u \in L_{p+d+e}H_{k+2d+2e}(Z',\mathbb{Q}).$$
\end{proposition}
\bp Let $X'=X/G_1$, $Y'=Y/G_2$, $Z'=Z/G_3$. Denote the three
quotient maps by $\pi_1: X\to X'$, $\pi_2: Y\to Y'$, $\pi_3: Z\to
Z'$. Define $\Gamma_1:=(\pi_1\times \pi_2)^*\Gamma'_1$ and
$\Gamma_2:=(\pi_2\times \pi_3)^*\Gamma'_2$. Consider the following
diagram (which looks like a prism with three square faces and two
triangular faces), our goal is to prove the bottom triangle commutes
on the level of $\mathbb{Q}$-coefficient Lawson homology groups.
 $$\xymatrix{X\ar[rr]^{\frac{\Gamma_2\circ\Gamma_1}{|G_2||G_3|}}\ar[rd]_(.7){\frac{\Gamma_1}{|G_2|}}\ar[dd]_{\pi_1}& & Z\ar[dd]^{\pi_3}\\
 & Y \ar[ru]_(.3){\frac{\Gamma_2}{|G_3|}}\ar[dd]^(.3){\pi_2}&\\
 X'\ar'[r]_(.7){\Gamma'_2\circ\Gamma'_1}[rr]\ar[rd]_{\Gamma'_1}& & Z'\\ & Y'\ar[ru]_{\Gamma'_2} &}$$
The upper triangle of the prism induces a commutative triangle in
$\mathbb{Q}$-coefficient Lawson homology by Proposition \ref{composition}.
The three squares also induce commutative squares in
$\mathbb{Q}$-coefficient Lawson homology. Indeed, for any $u\in
L_pH_k(X)$,
$$\begin{array}{lll}
\Gamma'_{1*}\pi_{1*}u&=p'_{2*}(p'^*_1\pi_{1*}u\cdot \Gamma'_1)\\
&=p'_{2*}\dfrac{(\pi_1\times\pi_2)_*}{|G_1||G_2|}\big{(}(\pi_1\times\pi_2)^*p'^*_1\pi_{1*}u\bullet (\pi_1\times\pi_2)^*\Gamma'_1\big{)}\\
&=\dfrac{\pi_{2*}p_{2*}}{|G_1||G_2|}\big{(}p_1^*\pi_1^*\pi_{1*}u\bullet \Gamma_1\big{)}\\
&=\dfrac{\pi_{2*}p_{2*}}{|G_1||G_2|}\big{(}p_1^*(\sum_{g\in
G_1}g_*u)\bullet \Gamma_1\big{)}
\end{array}
$$ where in the third equality we use the fact
$(\pi_1\times\pi_2)^*p_1'^*=p_1^*\pi_1^*$, which is valid even on
the level of cycles therefore valid on the level of Lawson homology.

Next, notice that for any $g\in G_1$, the identity
$p_1^*g_*=(g\times 1)_*p_1^*$ is valid on the level of cycles
therefore valid for Lawson homology. Moreover, $\Gamma_1$ is
invariant under the action of the group $(G_1\times 1)$. Therefore
by projection formula
$$p_{2*}(p_1^*g_*u\bullet\Gamma_1)=p_{2*}((g\times 1)_*p_1^*u\bullet\Gamma_1)=p_{2*}(g\times 1)_*(p_1^*u\bullet\Gamma_1)=p_{2*}(p_1^*u\bullet\Gamma_1).$$
Continue the above calculation of $\Gamma'_{1*}\pi_{1*}u$:
$$
\Gamma'_{1*}\pi_{1*}u
=\dfrac{\pi_{2*}}{|G_1||G_2|}\big{(}|G_1|p_{2*}p_1^*u\bullet
\Gamma_1\big{)}=\pi_{2*}\big{(}\frac{\Gamma_1}{|G_2|}\big{)}_*u.
$$
Thus the left square commutes. The commutativity of the other two
squares are similar, while in the proof we need fact that
$$\frac{\Gamma_2\circ\Gamma_1}{|G_2||G_3|}=(\pi_1\times \pi_3)^*
\big{(}\frac{\Gamma'_2\circ\Gamma'_1}{|G_3|}\big{)}.$$

Finally, since four of the five sides of the above prism induce
commutative diagrams and $\pi_{1*}: L_pH_k(X,\mathbb{Q})\to
L_pH_k(X',\mathbb{Q})$ is surjective, the triangle at the bottom
must commute, i.e.,
$$(\Gamma'_2\circ\Gamma'_1)_*=\Gamma'_{2*}\Gamma'_{1*}.$$
  \qe

In the similar situation as the last subsection, by the duality
isomorphism $\mathcal{D}'$, we have a corresponding result to
Proposition \ref{composition2} for the morphic cohomology follows
from Proposition \ref{quotient functorial}.

\begin{proposition}\label{quotient functorial2}
Let $X', Y', Z'$ be finite quotient varieties, $\Gamma'_1\in
  Corr^d(X',Y')$ and $\Gamma'_2\in Corr^e(Y', Z')$. Then for any $u\in L^qH^l(X',\mathbb{Q})$,
  $$(\Gamma'_2\circ\Gamma'_1)_*u=\Gamma'_{2*}\Gamma'_{1*}u \in L^{q+d+e}H^{l+2d+2e}(Z',\mathbb{Q}).$$
\end{proposition}\qe

\subsection{Motive, Lawson homology and morphic cohomology} In this subsection, we
explain that the morphic cohomology gives a covariant functor from
the category of Chow motives to the category of bi-graded abelian
groups. Analogously, the $\mathbb{Q}$-coefficient morphic cohomology
gives a covariant functor from the category of Chow motives for
finite quotient varieties to the category of bi-graded
$\mathbb{Q}$-vector spaces.


\medskip

We have recalled the definition of Chow motives in
\S\ref{sec:introduction} (Introduction). The theory of Chow motives
can be extended to $CH\mathcal{M'}$, the Chow motives of the
category of quotient varieties of smooth projective varieties by
finite groups (\cite{dBV}). To be more precise, let $\mathcal{V}'$
be the category of (not necessarily connected) varieties of the type
$X/G$ with $X\in Ob \mathcal{V}$ with an action of a finite group
$G$. The objects of $CH\mathcal{M}'$ are the same as the objects of
$\mathcal{V}'$, and the morphisms are defined similarly as in
$CH\mathcal{M}$. We again have a contravariant functor $h:
\mathcal{V}'\to CH\mathcal{M}'$.

\begin{theorem}\label{main1}
Given any $q, l\in \mathbb{Z}$, the morphic cohomology $L^qH^l$
defines a covariant functor from the category $\mathcal{CHM}$ to the
category of abelian groups as follows:
$$L^qH^l(X,\mathbf{p},r):=\mathbf{p}_*(L^{q+r}H^{l+2r}(X)) \subseteq
L^{q+r}H^{l+2r}(X).$$ Given a morphism $\Gamma\in
Hom_{CH\mathcal{M}}\big{(}(X,\mathbf{p},r),
(Y,\mathbf{q},s)\big{)}$, the morphism
$$L^qH^l(\Gamma): L^qH^l(X,\mathbf{p},r)\to L^qH^l(Y,\mathbf{q},s)$$ is defined as the
restriction of the map
$$\Gamma_*: L^{q+r}H^{l+2r}(X)\to L^{q+s}H^{l+2s}(Y).$$
\end{theorem}

\bp First, we need to show that $L^qH^l(\Gamma)$ is well defined,
i.e. the following diagram commutes,
$$\xymatrix{L^{q+r}H^{l+2r}(X)\ar[r]^{\Gamma_*}\ar[d]^{\mathbf{p}_*} & L^{q+s}H^{l+2s}(Y)\ar[d]^{\mathbf{q}_*}\\
L^{q+r}H^{l+2r}(X)\ar[r]^{\Gamma_*} & L^{q+s}H^{l+2s}(Y)}$$ It
commutes because of $\mathbf{q}\circ \Gamma=\Gamma\circ \mathbf{p}$
and Proposition \ref{composition2}.

Then we need to verify the functoriality of $L^qH^l$. This again
follows from Proposition \ref{composition2}.\qe

\begin{remark}
For our purpose, the category of Chow motives  $CH\mathcal{M}$ can
be replaced by the category of algebraic motives $CH_A\mathcal{M}$
whose objects are the same as $CH\mathcal{M}$, while morphisms are
defined to be
$$Hom_{CH_A\mathcal{M}}\big{(}(X,\mathbf{p},r), (Y,\mathbf{q},s)\big{)}:=\mathbf{q}\circ Corr^{s-r}_{alg}(X, Y)\circ \mathbf{p}$$
where $Corr^{s-r}_{alg}(X, Y)=Corr^{s-r}(X, Y)/\{{\rm algebraic ~ equivalence}\}$.
\end{remark}

\medskip

Similarly, for the Chow motives of finite quotient varieties we
have:

\begin{theorem}\label{main2}
  Given any $q, l\in \mathbb{Z}$, the $\mathbb{Q}$-coefficient morphic cohomology $L^qH^l(-,\mathbb{Q})$
defines a covariant functor from the category $CH\mathcal{M}'$ to
the category of abelian groups as follows:
$$L^qH^l((X',\mathbf{p},r),\mathbb{Q}):=\mathbf{p}_*(L^{q+r}H^{l+2r}(X',\mathbb{Q})) \subseteq
L^{q+r}H^{l+2r}(X',\mathbb{Q}).$$ Given a morphism $\Gamma\in
Hom_{CH\mathcal{M}'}\big{(}(X',\mathbf{p},r),
(Y',\mathbf{q},s)\big{)}$, the morphism
$$L^qH^l(\Gamma,\mathbb{Q}): L^qH^l((X',\mathbf{p},r),\mathbb{Q})\to L^qH^l((Y',\mathbf{q},s),\mathbb{Q})$$ is defined as the
restriction of map
$$\Gamma_*: L^{q+r}H^{l+2r}(X',\mathbb{Q})\to L^{q+s}H^{l+2s}(Y',\mathbb{Q}).$$
\end{theorem}
\bp Same as the proof of Theorem \ref{main1}. Proposition
\ref{quotient functorial2} implies that $L^qH^l(-,\mathbb{Q})$ is
well-defined and functorial.\qe

\medskip

There are corresponding versions of Theorem \ref{main1} and
\ref{main2} for Lawson homology, with almost the same proof hence we
skip it and only give the statement:

\begin{theorem}\label{Th4.3} For any $p, k\in \mathbb{Z}$,

(i) the Lawson homology $L_pH_k$ defines a contravariant functor
from the category $CH \mathcal{M}$ to the category of abelian groups
as follows:
$$L_pH_k(X,\mathbf{p},r):=\mathbf{p}_*(L_{p+r}H_{k+2r}(X)) \subseteq
L_{p+r}H_{k+2r}(X).$$   Given a morphism $\Gamma\in
Hom_{CH\mathcal{M}}\big{(}(X,\mathbf{p},r),
(Y,\mathbf{q},s)\big{)}$, the morphism
$$L_pH_k(\Gamma): L_pH_k(Y,\mathbf{q},s)\to L_pH_k(X,\mathbf{p},r)$$ is the
restriction of map $(^t\Gamma)_*: L_{p+s}H_{k+2s}(Y)\to
L_{p+r}H_{k+2r}(X).$

(ii) the $\mathbb{Q}$-coefficient Lawson cohomology
$L_pH_k(-,\mathbb{Q})$ defines a contravariant functor from the
category $CH \mathcal{M}'$ to the category of $\mathbb{Q}$-vector
spaces as follows:
$$L_pH_k((X',\mathbf{p},r),\mathbb{Q}):=\mathbf{p}_*(L_{p+r}H_{k+2r}(X',\mathbb{Q})) \subseteq
L_{p+r}H_{k+2r}(X',\mathbb{Q}).$$   Given a morphism $\Gamma\in
Hom_{CH\mathcal{M}'}\big{(}(X',\mathbf{p},r),
(Y',\mathbf{q},s)\big{)}$, the morphism
$$L_pH_k(\Gamma,\mathbb{Q}): L_pH_k((Y',\mathbf{q},s),\mathbb{Q})\to L_pH_k((X',\mathbf{p},r),\mathbb{Q})$$ is the
restriction of map $(^t\Gamma)_*: L_{p+s}H_{k+2s}(Y',\mathbb{Q})\to
L_{p+r}H_{k+2r}(X',\mathbb{Q}).$

\end{theorem}

\section{Applications}
\subsection{Projective bundles, blow-ups, and cell-decomposition}

As application of the connection between Lawson homology and the morphic cohomology,
we reobtain formulas for projective bundles, blow-ups, and smooth varieties admitting a
cell-decomposition. However, we require varieties to be smooth in these cases.

We start from the well known motivic decompositions for a projective
bundle and for a blow-up.  Let $\mathbb{P}$ be a projective bundle
over a smooth projective variety $X$ with fiber $\mathbb{P}^n$. The
following motivic decomposition is proved in \cite{Manin},
$$h(\mathbb{P})\simeq h(X)\oplus h(X)(1)\oplus\cdots\oplus h(X)(n).$$ then
by Theorem \ref{main1} and Theorem \ref{Th4.3} (recall that $(X,
\rm{id}_X, r)=h(X)(-r)$), we have the following result proved by
Friedlander and Gabber:
\begin{corollary}[Projective Bundle Theorem,
\cite{FG}]\label{projbundlethm} Let $\mathbb{P}$ be a projective
bundle over a smooth projective variety $X$ with fiber
$\mathbb{P}^n$. Then the following decompositions hold for morphic
cohomology and Lawson homology:
$$L^qH^l(\mathbb{P})\simeq L^qH^l(X)\oplus L^{q-1}H^{l-2}(X)\oplus\cdots\oplus L^{q-n}H^{l-2n}(X), \quad\forall q,l\in\mathbb{Z}.$$
\begin{equation}\label{eq4.1}  L_pH_k(\mathbb{P})\simeq
L_pH_k(X)\oplus L_{p-1}H_{k-2}(X)\oplus\cdots\oplus
L_{p-n}H_{k-2n}(X), \quad\forall p,k\in\mathbb{Z}.
\end{equation}

\end{corollary}

Let $X$ be a smooth projective variety and $j_0:V\hookrightarrow X$
a smooth subvariety of codimension $n\geq 2$. Let $\widetilde{X}$ be
the blowup of $X$ along $V$. Because of Theorem \ref{main1}, Theorem
\ref{Th4.3} and the motivic decomposition (cf. \cite{Manin})
$$h(\widetilde{X})\simeq h(X)\oplus h(V)(1)\oplus\cdots\oplus h(V)(n-1),$$
 we get the blowup formula for the morphic cohomology and Lawson homology:

\begin{corollary}[\cite{author1}]\label{cor4.2}
 Let $\widetilde{X}$ be the blow-up of a smooth projective variety
$X$ along a smooth subvariety $V$ of codimension $n$. Then
$$L^qH^l(\widetilde{X})\simeq L^qH^l(X)\oplus \bigoplus_{i=1}^{n-1}
L^{q-i}H^{l-2i}(V),$$
$$L_pH_k(\widetilde{X})\simeq L_pH_k(X)\oplus \bigoplus_{i=1}^{n-1}
L_{p-i}H_{k-2i}(V).$$
\end{corollary}

More generally, recall the following result proved by N. A. Karpenko
in \cite{Karpenko}:
\begin{theorem}[Karpenko] \label{celldecom}
Let $X$ be a smooth projective variety. Assume $X$  admits a filtration by closed subvarieties
$\emptyset=X_{-1}\subset X_{0}\subset\dots\subset X_{n}=X$
such that there exist flat morphisms $f_{i}: X_{i}-X_{i-1}\rightarrow Y_{i}$,
of relative dimension $m_{i}$ over smooth projective varieties $Y_{i}$ ($1\leq i\leq n$),
such that the fiber of every $f_{i}$ over every point $y$
of $Y_{i}$ is isomorphic to the affine space
$\C^{m_{i}}$. Then there exists an isomorphism
in $CH \mathcal{M}$
\begin{equation}\label{eq3}
 h(X)\simeq\bigoplus_{i=1}^{n} h(Y_{i})(m_{i}).
\end{equation}
 \end{theorem}

We immediately get the following
\begin{corollary} \label{cor4.3}
Using the notations in Theorem \ref{celldecom}, we have
\begin{equation}\label{eq4.2}
L^qH^l(X)\simeq\bigoplus_{i=1}^{n} L^{q-m_i}H^{l-2m_i}(Y_{i}),
\end{equation}
\begin{equation}\label{eq4.3}
  L_pH_k(X)\simeq\bigoplus_{i=1}^{n} L_{p-m_i}H_{k-2m_i}(Y_{i}).
\end{equation}
In particular, the isomorphism (\ref{eq4.3}) can be used to compute the
Lawson homology for  Grassmann bundles of projective vector bundles.
\end{corollary}

\bp
 Note that
\begin{equation}
\begin{array}{lll}
 Id_{X}&\in
 Hom_{CH\mathcal{M}}\big{(}(X,id,0), \oplus^{n}_{r=0}(Y_i,id,-m_i)\big{)}  \\
&=\bigoplus^{n}_{i=0}Hom_{CH\mathcal{M}}\big{(}(X,id,0), (Y_i,id,-m_i)\big{)} \\
&=\bigoplus^{n}_{i=0}Corr^{-m_i}(X,Y_i), \\
\end{array}
\end{equation}
hence $Id_{X}=\oplus^{n}_{i=0}\Gamma_i$, where
$\Gamma_i\in Corr^{-m_i}(X,Y_i)
$.
By Theorem
\ref{main1} and \ref{celldecom}, we have Equation (\ref{eq4.2}).

Similarly, by Theorem \ref{Th4.3} and \ref{celldecom}, we obtain Equation (\ref{eq4.3}).
\qe

\begin{remark}
It was pointed out to us by Friedlander that the decomposition of motives in Equation (\ref{eq3}) implies
the decomposition of any oriented cohomology theory(cf. \cite{Nenashev-Zainoulline}). Friedlander and Walker
showed that the Lawson homology and morphic cohomology are such theories for varieties over
$\R$ and $\C$ (cf. \cite{FW2} and references therein by
the same authors). Hence, Corollary \ref{cor4.3} was implied from those,
although explicit formula was not written down.
\end{remark}

\begin{remark}
By the Friedlander-Lawson duality (cf. Proposition \ref{negative}), the Equation
(\ref{eq4.2}) is equivalent to  the following formula in terms of Lawson homology groups:
\begin{equation}\label{eq4.4}
L_pH_k(X)\simeq\bigoplus_{i=1}^{n} L_{p-d_i}H_{k-2d_i}(Y_{i}).
\end{equation}
By comparing Equation (\ref{eq4.3}) and (\ref{eq4.4}), we obtained visible obstructions
for a collection of pairs $\{(Y_i,m_i)\}_{i=1}^{n}$, where $Y_i$ is a smooth projective
variety and $m_i$ is a positive integer for each $1\leq i\leq n$, to be the decomposition
of a smooth projective variety $X$ in the sense of Theorem \ref{celldecom}.

\end{remark}


\subsection{Hilbert scheme of points on a
surface}\label{sec:Hilbert}

{\hskip .2 in} \label{sec:3} In this section, we will compute the
rational coefficient morphic cohomology and Lawson homology for
Hilbert scheme of points on a smooth complex projective surface.

Let $X$ be a smooth projective surface, let $X^{(n)}$ be its $n$-th
symmetric product, let $X^{[n]}$ be the Hilbert scheme of
0-dimensional subschemes of $X$ of length $n$, and let
$\pi:X^{[n]}\rightarrow X^{(n)}$ be the Hilbert-Chow morphism. It is
well-known that $X^{[n]}$ is nonsingular. We denote by
$\mathfrak{P}(n)$ the set of partitions of $n$ and $p(n)$ its
cardinality. For any $\nu\in \mathfrak{P}(n)$, we denote by $l(\nu)$
its length, and define $X_\nu^{(n)}$ to be the locally closed subset
of points in $X^{(n)}$ of the type
${\nu}_1x_1+\cdots+{\nu}_{l(\nu)}x_{l(\nu)}$, with $x_h\in X$ and
$x_i\neq x_j$ for every $i\neq j$. Define $X_\nu^{[n]}$ to be the
reduced scheme $(\pi^{-1}(X_\nu^{(n)}))_{red}$. Let
$\overline{X}_\nu^{[n]}$ be the closure of the stratum $X_\nu^{[n]}$
in $X^{[n]}$ and let $\overline{X}_\nu^{(n)}$ be the closure of
$X_\nu^{(n)}$ in $X^{(n)}$. It can be proved that
$\overline{X}_\nu^{[n]}=\pi^{-1}(\overline{X}_{\nu}^{(n)})$. If
$\nu=1^{a_1}\cdots n^{a_n}$, then the finite group
$\Sigma_{\nu}:=\Sigma_{a_1}\times\cdots\times \Sigma_{a_n}$ acts
naturally on $X^{l({\nu})}$, where $\Sigma_{a_i}$ are the symmetric
groups. The quotient $X^{\nu}$ is isomorphic to
$X^{(a_1)}\times\cdots\times X^{(a_n)}$. Use the notations in
\cite{dCM}, we denote $X^{l({\nu})}$ by $X^{\nu}$. The natural
$\Sigma_{{\nu}}$-invariant map ${\nu}:X^{\nu}\to X^{(n)}$ has image
$\overline{X}_{\nu}^{(n)}$. Hence it descends to a map
${\nu}:X^{({\nu})}\to X^{(n)}$ which we denote by the same symbol.
By using these notations, the correspondences $\Gamma^{\nu}$ and
$\widehat{\Gamma}^{\nu}$ are defined as follows:
$$
  \Gamma^{\nu}:=\{(x_1,\cdots,x_{l({\nu})},\J)\in X^{\nu}\times X^{[n]}:
\pi(\J)={\nu}_1x_1+\cdots {\nu}_{l({\nu})}x_{l({\nu})})
\}\cong \big(X^{\nu}\times_{X^{(n)}} X^{[n]}\big)_{red}.
  $$
and
$$
\widehat{\Gamma}^{\nu}:=\Gamma^{\nu}/\Sigma_{\nu}
$$
since the
correspondence $\Gamma^{\nu}$ is invariant under the action of
$\Sigma_{\nu}$ on the first factor of the product.

Set
$\widehat{{\mathcal{X}}}=\coprod_{\nu\in\mathfrak{P}(n)}X^{(\nu)}$
and
$\widehat{{\Gamma}}=\coprod_{\nu\in\mathfrak{P}(n)}\widehat{\Gamma}^{\nu}$.
Define integer $m_\nu:=(-1)^{n-l(\nu)}\prod_{j=1}^{l{(\nu)}}\nu_j$.
Let
$$
\widehat{\Gamma}':=\bigoplus_{\nu\in
\mathfrak{P}(n)}\frac{^t\widehat{\Gamma}^{\nu}}{m_\nu}
$$
where
$^t\widehat{\Gamma}^{\nu}$ is the transposed correspondence of
$\widehat{\Gamma}^{\nu}$. It is proved by de Cataldo and Migliorini
that
\begin{theorem}[\cite{dCM}]\label{dCM}
The correspondence $\widehat{{\Gamma}}$ gives an isomorphism of
Chow motives (in the category of finite quotient varieties):
 $$
 \widehat{{\Gamma}}=\bigoplus_{\nu\in \mathfrak{P}(n)}\widehat{\Gamma}^{\nu}:
 \bigoplus_{\nu\in \mathfrak{P}(n)}(X^{(\nu)}, \Delta_{X^{(\nu)}})(n-l(\nu))\longrightarrow (X^{[n]}, \Delta_{X^{[n]}}),
 $$
with the inverse correspondence given by $\widehat{\Gamma}'$.
\end{theorem}

As a consequence, we obtain the following theorem:

\begin{theorem}\label{main3}
{Let $X$ be a smooth complex projective surface. The natural map of
morphic cohomology groups
\begin{equation}\label{eq main1}\widehat{{\Gamma}}_*=\bigoplus_{\nu\in
\mathfrak{P}(n)}\widehat{\Gamma}^\nu_*: \bigoplus_{\nu\in
\mathfrak{P}(n)}L^{q-n+l(\nu)}H^{l-2n+2l(\nu)}(X^{(\nu)},
\mathbb{Q})\longrightarrow L^qH^l(X^{[n]},\mathbb{Q})
\end{equation}
is an isomorphism for any integers $0\leq l\leq 2q$. (For other $l,
q$ both sides becomes zero, therefore the isomorphism holds
trivially).

Similarly, we have an isomorphism of Lawson homology groups for all
$p,k\in \mathbb{Z}$:
\begin{equation}\label{eq main2}\bigoplus_{\nu\in
\mathfrak{P}(n)}L_{p-n+l(\nu)}H_{k-2n+2l(\nu)}(X^{(\nu)},
\mathbb{Q})\longrightarrow L_pH_k(X^{[n]},\mathbb{Q})
\end{equation}
}
\end{theorem}
\bp Note that the $\dim \widehat{{\Gamma}}^\nu=\dim
{{\Gamma}}^\nu=n+l(\nu)$. By the functoriality proved in Theorem
\ref{main2}, we have
$\widehat{{\Gamma}}_*\widehat{{\Gamma}}_*'=(\widehat{{\Gamma}}\widehat{{\Gamma}}')_*={\rm
id_*}$,
$\widehat{{\Gamma}}_*'\widehat{{\Gamma}}_*=(\widehat{{\Gamma}}'\widehat{{\Gamma}})_*={\rm
id_*}$. Therefore $\widehat{{\Gamma}}_*$ gives an isomorphism
between morphic cohomology groups by Theorem \ref{dCM}, this proves
(\ref{eq main1}). The isomorphism (\ref{eq main2}) is obtained by
applying Theorem \ref{Th4.3} to Theorem \ref{dCM}.

\qe
\begin{remark}
A direct proof of the above theorem (without using the language of
motives)
 can be obtained from the method in \cite{dCM} to compute the Chow group of $X^{[n]}$.
\end{remark}

\begin{remark}
 The duality isomorphism between morphic cohomology and Lawson
homology can be used to prove the equivalence of the two isomorphisms in
Theorem \ref{main3}.
\end{remark}

The above result can be applied to Friedlander-Walker semi-topological $K$-theory (cf. \cite{FW2} and references therein). Notice that, by
\cite{Fulton} Corollary 18.3.2, de Cataldo and Migliorini (in
\cite{dCM} Theorem 5.4.1) give a decomposition of the rational
coefficient Grothendieck group $K_0(X^{[n]})_{\mathbb
Q}(:=K_0(X^{[n]})\otimes {\mathbb  Q})$:
\begin{equation}\label{eq5.1}
 \widehat{{\Gamma}}_*=\bigoplus_{\nu\in \mathfrak{P}(n)}\widehat{\Gamma}^{\nu}_*:
\bigoplus_{\nu\in\mathfrak{P}(n)}{K}_{0}(X^{(\nu)})_{\mathbb{Q}}\stackrel{\cong}{\rightarrow}
{ K}_0(X^{[n]})_{\mathbb{Q}}
\end{equation}
They asked if similar statements hold for higher $K$-theory.


We do not have an answer for this question. Instead, we give an
answer to a similar question for the semi-topological $K$-theory.

The following is an immediate consequence of Theorem \ref{main3}. It gives a decomposition
to $K_*^{sst}(X^{[n]})_{_{\mathbb{Q}}}$ in terms of rational Lawson homology groups.

\begin{corollary}
There is a natural isomorphism of the  semi-topological $K$-theory groups with
rational coefficients
 $$
{K}_p^{sst}(X^{[n]})_{\mathbb{Q}}  \stackrel{\cong}{\longrightarrow} \bigoplus_{\nu\in\mathfrak{P}(n)}
\bigoplus_j{L}_{j}H_{2j+p-n+l(\nu)}(X^{(\nu)})_{\mathbb{Q}}.
 $$
\end{corollary}
\bp It follows from Theorem 4.7 in \cite{FW} and Theorem \ref{main3}.\qe

\begin{remark} We expect the following isomorphism holds
 $$
\bigoplus_{\nu\in\mathfrak{P}(n)}{K}_{p-n+l(\nu)}^{sst}(X^{(\nu)})_{\mathbb{Q}}\rightarrow {
K}_p^{sst}(X^{[n]})_{\mathbb{Q}}
 $$
 for a smooth projective surface $X$.
 However we need a similar isomorphism as that in Theorem 4.7 in \cite{FW} for finite quotient varieties, or
 more specifically, for $X^{(\nu)}$.
\end{remark}

\section{Further consequences}

\subsection{Birational invariants defined by Lawson homology using correspondence}   The action of correspondences on Lawson
homology gives another proof of the following theorem, which is
originally discovered in \cite{author1} by using diagram chases and
blow up formula for Lawson homology together with the Weak
Factorization Theorem.

\begin{theorem}[\cite{author1}]\label{birat-inv}
If $X$ and $Y$ are birationally equivalent smooth projective complex
varieties of dimension $n$, then $L_1H_k(X)_{hom}\cong
L_1H_k(Y)_{hom}$ for $k\geq 2$ and $L_{n-2}H_k(X)_{hom}\cong
L_{n-2}H_k(Y)_{hom}$ for $k\geq 2(n-2)$.
\end{theorem}
\bp We prove only the isomorphism for $L_1H_k(-)_{hom}$. The proof
of the isomorphism for $L_{n-2}H_k(-)_{hom}$ is similar.

Let $\Gamma$ be the closure of the graph of a birational map $f:
X\dashrightarrow Y$. Note that $\Gamma^t\circ\Gamma$ is the sum of
identity correspondence $\Delta_X$ and correspondences $\gamma_i$'s
whose projections are contained in proper subvarieties of  $X$ (cf.
Example 16.1.11 in \cite{Fulton}). Indeed, let $U\subseteq X$,
$U'\subseteq Y$ such that $f$ restricts to $U$ is an isomorphism to
$U'$. Let $D_X=X\setminus U$ and $D_Y=Y\setminus U'$. It is easy to
see that $(\Gamma^t\circ\Gamma-\Delta_X)$ can be chosen to be
supported in $D_X\times D_X$.

By a result of Peters (Lemma 12 in \cite{Peters}, cf. Remark
\ref{peters} below), for any $u\in L_1H_k(X)_{hom}$ and $\gamma\in
\mathcal{C}_n(X\times Y)$ such that $p_1(\gamma)$ is a proper
subvariety in $X$, we have $\gamma_*(u)=0\in L_1H_k(Y)_{hom}$.
Therefore, $\Gamma^t_*\Gamma_*=(\Gamma^t\Gamma)_*=(\Delta_X)_*=id$
on $L_1H_k(X)_{hom}$. Symmetrically, $\Gamma_*\Gamma^t_*=id$ on
$L_1H_k(Y)_{hom}$. Therefore $\Gamma_*$ induces an isomorphism
$$L_1H_k(X)_{hom}\cong L_1H_k(Y)_{hom}.$$\qe

\begin{remark}\label{peters}
Lemma 12 in \cite{Peters} asserts that: assume $X$ and $Y$ are
smooth projective varieties and $\alpha\subset X\times Y$ is an
irreducible cycle of dimension $\dim X=n$, supported in $V\times W$
where $\dim V=v$ and $\dim W=w$. Then $\alpha_*=0$ if $m<n-v$ or if
$m>w$. Moreover, $\alpha_*=0$ on $L_{n-v}H_*(X)_{hom}$ and on
$L_wH_*(X)_{hom}$.

The statement that ``$\alpha_*=0$ if $m<n-v$'' is not correct, since
$L_pH_k(-)$ is not necessarily zero. But other statements are still
valid which we explain here:
Let $\tilde{V}\to V$ and $\tilde{W}\to
W$ be resolutions of singularities. Let $i:\tilde{V}\to X$ and $j:
\tilde{W}\to Y$ be the natural morphisms, $\tilde{\alpha}$ be the
proper transform of $\alpha$, it can be checked that  the following
diagram commutes (cf. the proof of Lemma 12 in \cite{Peters}):
$$\xymatrix{L_{m-n+v+w}H_{l+2(v+w-n)}(\tilde{V}\times\tilde{W})_{hom}\ar[r]^-{\tilde{\alpha}_*}& L_mH_l(\tilde{V}\times\tilde{W})_{hom}\ar[d]_{(p_2)_*}\\
    L_{m-n+v}H_{l+2(v-n)}(\tilde{V})_{hom}\ar[u]^{p_1^*}& L_mH_l(\tilde{W})_{hom}\ar[d]_{j_*} \\
    L_mH_l(X)_{hom}\ar[u]^{i^*}\ar[r]^{\alpha_*}&  L_mH_l(Y)_{hom}\\}$$
As a consequence of this commutative diagram, the conclusion that
$\alpha_*=0$ on $L_{n-v}H_*(X)_{hom}$ is valid by noticing that
$L_pH_k(-)_{hom}=0$ for $p\leq 0$.
\end{remark}

\subsection{Unirational threefolds and more}

In this subsection we describe the Lawson homology for unirational
threefolds and fourfolds, and more general the relation between the
Lawson homologies of two varieties $X$ and $Y$.

First of all, we make a remark the motive of a curve. Given a smooth
projective curve $C$ and a point $e\in C$, we put
$\textbf{p}_0=e\times C$ and $\textbf{p}_2=C\times e$, then take
$\textbf{p}_1=\Delta_C-\textbf{p}_0-\textbf{p}_2$ where $\Delta_C$
is the diagonal in $C\times C$. Then we have $h(C)=h(pt)\oplus
\L\oplus C^+$, where $\L=h(pt)(1)$ is the Lefschetz motive and
$C^+=(C,id-\textbf{p}_0-\textbf{p}_2)$. It is known that the natural
map $L_pH_k(C)\to H_k(C)$, namely the cycle map, is an isomorphism.
It is also easy to show that the cycle map commutes with the map
$\Gamma_*$ induced from any correspondence $\Gamma$. Therefore
$L_pH_k(C^+)\to H_k(C^+)$ is also an isomorphism.

In \S 11 of Manin's paper \cite{Manin}, he gives a motivic
decomposition of a unirational threefold $X$, namely
$$  h(X)=h(pt)\oplus a\L\oplus U\otimes\L\oplus a\L^2\oplus\L^3,
$$
where $U$ is a direct summand of a motive of the form $\oplus
Y_i^+$, the $Y_i$ being curves. By the argument in the previous
paragraph,  $L_pH_k(U,\mathbb{Q})\cong H_k(U,\mathbb{Q})$. Then by
Theorem \ref{Th4.3} (i), a motive decomposition implies the
decomposition of rational Lawson homology as well as it is
well-known for the rational singular homology, we obtain the
following:
\begin{proposition} Let $X$ be a three
dimensional smooth projective unirational
 variety over $\C$. Then the rational Lawson homology group
 is isomorphic to the corresponding rational singular homology groups,i.e.,
$$L_pH_k(X,\mathbb{Q})\cong H_k(X,\mathbb{Q})
$$for any $p$ and $k$.
\end{proposition}

By a similar argument, we obtain the following result for a
unirational fourfold:
\begin{proposition}
Let $X$ be a unirational smooth complex projective variety of
dimension four. Then the relation of  rational Lawson homology and
rational singular cohomology is given as follows:
$$
 \left\{\begin{array}{ll}
 L_pH_k(X,\mathbb{Q}) \cong H_k(X,\mathbb{Q}), & \hbox{if $(p,k) \neq (2,4)$;} \\
 L_pH_k(X,\mathbb{Q}) \hookrightarrow H_k(X,\mathbb{Q}) ~\hbox{is injective}, & \hbox{if $(p,k)=(2,4)$.}
 \end{array}
 \right.
$$
\end{proposition}

\begin{remark}
These results can be obtained by the method on the decomposition
of diagonals, used by C. Peters in
\cite{Peters}, since any unirational variety has small Chow groups
for zero cycles. Later, M. Voineagu refines Peters' results, which
consider only $\mathbb{Q}$-coefficient Lawson homology groups, to
$\mathbb{Z}$-coefficient in many cases \cite{Voineagu}.
\end{remark}

The above propositions can be generalized to a generically finite
rational map as follows:
\begin{proposition}\label{finite-map}
If $f: X\dashrightarrow Y$ be a generically finite rational map
between smooth projective varieties of dimension $n$. Then
\begin{equation}\label{eq9}
\dim_{\mathbb{Q}}L_1H_k(X,\mathbb{Q})_{hom}\geq
\dim_{\mathbb{Q}}L_1H_k(Y,\mathbb{Q})_{hom}
\end{equation}
for $k\geq 2$ and
\begin{equation}\label{eq10}
\dim_{\mathbb{Q}}L_{n-2}H_k(X,\mathbb{Q})_{hom}\geq
\dim_{\mathbb{Q}}L_{n-2}H_k(Y,\mathbb{Q})_{hom}
\end{equation}
for $k\geq 2(n-2)$. (In case that the right hand side of the
inequality has infinite dimension, the left hand side must also be
infinite dimensional.)
\end{proposition}

\bp Note that $f: X\dashrightarrow Y$ is  a rational map of degree
$d>0$. By Hironaka's theorem on the resolution of singularities of
mappings,  there is a commutative diagram of the form
$$\xymatrix{\widetilde{X}\ar[d]^{\sigma}\ar[rd]^F&\\
X\ar@{.>}[r]^f &Y }
$$
where $F$ is a morphism of degree $d$, and  $\sigma$ is a birational
morphism which factors into  the composition of blow ups on smooth
subvarieties of codimension at least $2$. It can be proved that
\begin{equation}\label{eq11}
\sigma_*:L_1H_k(\widetilde{X},\mathbb{Q})_{hom}\to
L_1H_k(X,\mathbb{Q})_{hom}
\end{equation}
(even $L_1H_k(\widetilde{X})_{hom}\to L_1H_k(X)_{hom}$) is an
isomorphism by reducing the result to one blow up as given in
\cite{author1}. The same is true for codimension two cycles.

It remains to prove that
$F_*:L_1H_k(\widetilde{X},\mathbb{Q})_{hom}\to
L_1H_k(Y,\mathbb{Q})_{hom}$ is surjective. Let $\Gamma_F\in
\Ch_n(\widetilde{X}\times Y)$ be the graph of $F$ and $\Gamma_F^t$ be
its transpose.  Since $F:\widetilde{X}\to Y$ is a morphism of finite
degree $d$ between smooth projective varieties, we have
$h(\widetilde{X})=h(Y)\oplus
(\widetilde{X},id_{\widetilde{X}}-\mathbf{p},0)$, where
$\mathbf{p}=\frac{1}{d}(\Gamma_F^t)\circ(\Gamma_F)$. Therefore we
have
$\mathbf{p}_*(L_pH_k(\widetilde{X},\mathbb{Q}))=L_pH_k(Y,\mathbb{Q})$
and $\mathbf{p}_*(H_k(\widetilde{X},\mathbb{Q}))=H_k(Y,\mathbb{Q})$.
These two equations imply
$\mathbf{p}_*(L_pH_k(\widetilde{X},\mathbb{Q})_{hom})=L_pH_k(Y,\mathbb{Q})_{hom}$
since pull-backs and push-forwards commute with the natural
transformation from the Lawson homology to the singular homology.
Therefore
\begin{equation}\label{eq12}
\dim_{\mathbb{Q}}L_pH_k(\widetilde{X},\mathbb{Q})_{hom}\geq
\dim_{\mathbb{Q}}L_pH_k(Y,\mathbb{Q})_{hom}
\end{equation}

From Equations (\ref{eq11}) and (\ref{eq12}), we get Equation
(\ref{eq9}). Similar for Equation (\ref{eq10}).

\qe

In particular, for a uniruled threefold $X$, (recall that a threefold $X$ is uniruled if
there is a generic finite map $f:S\times \P^1\dashrightarrow X$ for some surface
$S$)
$L_pH_k(X,\Q)\cong H_k(X,\Q)$ if $(p,k)\neq (1,2)$ or $(2,4)$ and
$L_pH_k(X,\Q)\hookrightarrow H_k(X,\Q)$ is injective if $(p,k)= (1,2)$ or $(2,4)$.

\begin{remark}\label{finite-morphism}
From the proof of Proposition \ref{finite-map}, we see that if we have a finite morphism
$f:X\to Y$ between smooth projective varieties, then
\begin{equation}\label{eq15}
\dim_{\mathbb{Q}}L_pH_k(X,\mathbb{Q})_{hom}\geq
\dim_{\mathbb{Q}}L_pH_k(Y,\mathbb{Q})_{hom}.
\end{equation}
\end{remark}

\subsection{Griffiths groups for the product of curves}

As the application of the above Proposition \ref{finite-map},
together results on Griffiths group on generic Jacobian of smooth projective curves,
we give examples of products of smooth curves carrying
nontrivial Griffiths groups.

\begin{proposition}\label{nonzeroGriff}
Let $C$ be generic smooth projective curve of genus $g\geq 3$ and let $X=C^g$ be the
$g$-copies of self products of $C$.
Then $\Griff_p(X)\otimes\Q$ are nontrivial for all $1\leq p\leq g-2$.
\end{proposition}
\bp
Let $C$ be a generic curve of genus $g\geq 3$. Firstly, note that the Jacobian
$J(C)$ of $C$ have a non-trivial Griffiths group $\Griff_p(J(C))\otimes \mathbb{Q}$ for $1\leq p\leq g-2$ (cf. \cite{Ceresa}).

Secondly, it is well known that there is a birational morphism from
the $g$-th symmetric product $C^{(g)}$ of $C$ to  $J(C)$, i.e., $\sigma: C^{(g)}\to J(C)$ is
a birational morphism. Therefore,
$$\dim_{\Q}\{\Griff_p(C^{(g)})\otimes\Q\}\geq \dim_{\Q}\{\Griff_p(J(C))\otimes \mathbb{Q}\}$$
by the proof to Equation (\ref{eq15}) in Remark \ref{finite-morphism}.
For the special cases $p=1$ or $g-2$,  $\Griff_1(C^{(g)})\cong \Griff_1(J(C))$ and
$\Griff_{g-2}(C^{(g)})\cong \Griff_{g-2}(J(C))$  also follows from Theorem \ref{birat-inv}.

Finally, since the natural projection $\pi:C^g\to C^{(g)}$
is of finite degree $g!$, we have
$$\dim_{\Q}\{\Griff_p(C^{g})\otimes\Q\}\geq \dim_{\Q} \{\Griff_p(C^{(g)})\otimes \Q\}$$
 from Remark \ref{finite-morphism}.

The combination these statements completes the proof of the proposition.
\qe

\begin{remark}
Since  all Griffiths groups for curves are zero and
$L_1H_2(X)_{hom}\cong \Griff_1(X)$ for a smooth projective variety $X$, we obtain that
the K\"{u}nneth type formula in general can \textbf{not}
hold for Griffiths groups and Lawson homology.
\end{remark}

\begin{remark}
It was constructed explicitly by B. Harris in \cite{B-Harris} for the Fermat curve $C$ of
degree $4$(hence $g(C)=3$) with $\Griff_1(J(C))\neq0$. In fact, $\Griff_1(J(C))\otimes \Q\neq 0$.
By the proof in Proposition  \ref{nonzeroGriff}, we get $\Griff_1(X)\otimes \Q\neq 0$ for $X=C^3$ the $3$ times self product
of $C$.
\end{remark}

\section{Appendix}
We would like to make a remark on birational morphisms.

\begin{proposition}
Let $f:X\to Y$ be a morphism between smooth complex projective
$n$-dimensional varieties. Suppose that $f_*:\Ch_p(X)\to \Ch_p(Y)$
is isomorphic for $p=n-1$ and $p=n$. Then $f$ is an isomorphism.
\end{proposition}
\bp The isomorphism of $f_*$ for $p=n$ implies that $f$ is a
birational morphism. Indeed, $\Ch_n(X)\cong\Ch_n(Y)\cong\mathbb{Z}$,
and $f_*$ is an multiplication by $d$ where $d$ is the degree of
$f$, hence $d=1$. Then we apply the fact that, if a birational
morphism $f$ is not an isomorphism, then there is an exceptional
subvariety $Z\subset X$, i.e. $\textrm{codim } Z=1$ and
$\textrm{codim } f(Z)\ge 2$. It is easy to show that $[Z]\neq 0\in
\Ch_{n-1}(X)$ but $f_*([Z])=0\in \Ch_{n-1}(Y)$, contradicts to the
fact that $f_*$ is an isomorphism.

\qe

\begin{remark}
The statement is amazing by comparing to the corresponding one between topological manifolds:
For a continuous map $F:M\to N$ between oriented topological manifolds,
even if $F_*: H_k(M,\Z)\to H_k(N,\Z)$ are isomorphisms for all $k$,
we don't know whether $F:M\to N$ is homeomorphic.

\end{remark}

\medskip

\noindent Wenchuan Hu, Department of Mathematics, MIT, Room 2-101,
77 Massachusetts Avenue Cambridge, MA 02139 \quad Email:
wenchuan@math.mit.edu

\medskip

\noindent Li Li, Department of Mathematics, Univ. of Illinois,
Urbana, IL 61801 \\ Email: llpku@math.uiuc.edu


\begin{thebibliography}{AAAA}

\bibitem[Ce]{Ceresa}G. Ceresa, {\sl $C$ is not algebraically equivalent to $C\sp{-}$ in its Jacobian.}
Ann. of Math. (2) 117 (1983), no. 2, 285--291.


\bibitem[CH]{CH} A.~Corti, M.~Hanamura, {\em Motivic decomposition and intersection Chow groups. I}, Duke
Math. J. 103 (2000), no. 3, 459--522.

\bibitem[dBV]{dBV}S.~del Ba\~{n}o, V.~Navarro, {\em On the motive of a
quotient variety}, Dedicated to the memory of Fernando Serrano.
Collect. Math. 49 (1998), no. 2-3, 203--226.

\bibitem[dCM]{dCM}
M.~de Cataldo and L.~Migliorini, {\em The Chow groups and the motive
of the Hilbert scheme of points on a surface.} J. Algebra 251
(2002), no. 2, 824--848.

\bibitem[Fo]{Fogarty} J. Fogarty,
{\sl Algebraic families on an algebraic surface.}
Amer. J. Math 90 1968 511--521.


\bibitem[FHW]{Friedlander-Haesemesyer-Walker}Friedlander, Eric M.(1-NW); Walker, Mark E.(1-NE)
Semi-topological $K$-theory. Handbook of $K$-theory. Vol. 1, 2,
877--924, Springer, Berlin, 2005.

\bibitem[F]{Friedlander1} E. Friedlander, {\sl Algebraic cycles, Chow
varieties, and Lawson homology.}  Compositio Math. 77 (1991), no. 1,
55--93.

\bibitem[FG]{FG} E. Friedlander and O. Gabber, {\sl Cycle
spaces and intersection theory. Topological methods in modern
mathematics.} (Stony Brook, NY, 1991), 325--370, Publish or Perish,
Houston, TX, 1993.


\bibitem[FL1]{FL1}E. Friedlander and B. Lawson {\sl A theory of algebraic cocycles.}
Ann. of Math. (2)  136  (1992),  no. 2, 361--428.

\bibitem[FL2]{FL2}E. Friedlander and B. Lawson {\sl Duality relating spaces of
algebraic cocycles and cycles.}  Topology  36  (1997),  no. 2,
533--565.


\bibitem[FW]{FW}E. Friedlander and M. Walker, {\sl
Rational isomorphisms between $K$-theories and cohomology theories.}
Invent. Math. 154 (2003), no. 1, 1--61.

\bibitem[FW2]{FW2} E. Friedlander and M. Walker,
{\sl Semi-topological $K$-theory.} Handbook of $K$-theory. Vol. 1, 2, 877--924, Springer, Berlin, 2005.

\bibitem[Fu]{Fulton}
W. Fulton, Intersection theory. Second edition, Springer-Verlag,
Berlin, 1998.

\bibitem[Hb]{B-Harris}
B. Harris,
{\em Homological versus algebraic equivalence in a Jacobian.}
Proc. Nat. Acad. Sci. U.S.A. 80 (1983), no. 4 i., 1157--1158.

\bibitem[Hj]{Harris} J. Harris, {\sl Algebraic geometry, A first course}, Graduate Texts in
Mathematics, 133. Springer-Verlag, New York, 1992.

\bibitem[Hu]{author1} W. Hu, {\sl Birational invariants defined by
Lawson homology.}\\  arXiv:math.AG/0511722.

\bibitem[K]{Karpenko}  N. A. Karpenko,
{\sl Cohomology of relative cellular spaces and of isotropic flag varieties.}
(Russian. Russian summary) Algebra i Analiz 12 (2000), no. 1, 3--69; translation in
St. Petersburg Math. J. 12 (2001), no. 1, 1--50



\bibitem[L1]{Lawson1}
B. Lawson, {\sl Algebraic cycles and homotopy theory.}, Ann. of
Math. {\bf 129}(1989), 253-291.

\bibitem[L2]{Lawson2}B. Lawson, {\sl Spaces of algebraic
cycles.} pp. 137-213 in Surveys in Differential Geometry, 1995
vol.2, International Press, 1995.

\bibitem[LF]{Lima}P. Lima-Filho, {\sl Lawson homology for quasiprojective
varieties.} Compositio Math.  84  (1992),  no. 1, 1--23.

\bibitem[LF2]{Lima2} P. Lima-Filho, {\sl The topological group structure of algebraic
cycles.} Duke Math. J. 75 (1994), no. 2, 467--491.

\bibitem[LF3]{Lima3}
Lima-Filho, Paulo
On the generalized cycle map. (English summary)
J. Differential Geom. 38 (1993), no. 1, 105--129.

\bibitem[M]{Manin}
J. I. ~Manin, {\em Correspondences, motifs and monoidal
transformations}, Math. USSR-Sbornik 6 (1968), 439--470.


\bibitem[NZ]{Nenashev-Zainoulline}
A. Nenashev and K. Zainoulline,
{\em Oriented cohomology and motivic decompositions of relative cellular spaces.} (English summary)
J. Pure Appl. Algebra 205 (2006), no. 2, 323--340.




\bibitem[Pe]{Peters} C. Peters, {\sl  Lawson homology for varieties with
small Chow groups and the induced filtration on the Griffiths
groups.}  Math. Z. 234 (2000), no. 2, 209--223.


\bibitem[Vo]{Voineagu} Mircea Voineagu,
{\sl Semi-topological K-theory for certain projective varieties.}
Preprint. arxiv.org/abs/math/0601008

\end{thebibliography}
\end{document}